\documentclass[11pt,a4paper]{article}
\usepackage{tikz}
\usetikzlibrary{patterns}
\usepackage{graphicx}
\usepackage{subfigure}
\usepackage{xcolor}
\usepackage{fullpage}

\usepackage{amsmath,amsthm,amsfonts}  
\newtheorem{theorem}{Theorem}[section]

\newtheorem{remark}[theorem]{Remark}

\newcommand{\Om}{\Omega}
\newcommand{\bR}{\mathbb{R}}
\newcommand{\Ga}{\Gamma}
\renewcommand{\div}{\textrm{div}}
\newcommand{\D}[2]{\frac{\partial #1}{\partial#2}}
\renewcommand{\d}{\diamond}
\newcommand{\ga}{\gamma} 
\newcommand{\Ld}{\Lambda}
\newcommand{\ld}{\lambda}
\newcommand{\wt}{\widetilde}
\newcommand{\ve}{\varepsilon}

\begin{document}

\title{The Factorization method for three dimensional Electrical Impedance Tomography}

\author{N. Chaulet\footnote{Department of Mathematics University College London,  London WC1E 6BT, UK}, 
S. Arridge\footnote{Centre for Medical Image Computing, University College London, London WC1E 6BT, UK},
 T. Betcke$^*$ 
 and D. Holder\footnote{Department of Medical Physics and Bioengineering, University College London,  London 
WC1E 6BT, UK}}
\date{}
\maketitle
\begin{abstract}
The use of the Factorization method for Electrical Impedance Tomography has been proved to be very promising for applications in the case where one wants to 
find inhomogeneous inclusions in a known background. In many situations, the inspected domain is three dimensional and is made of various materials. In this 
case, the 
main challenge in applying the Factorization method consists in computing the Neumann Green's function of the background medium. We explain how we solve 
this difficulty and demonstrate the 
capability of the Factorization method to locate inclusions in realistic inhomogeneous three dimensional background media from  simulated data obtained 
by 
solving the so-called complete electrode model. 
We also perform a numerical study of the stability of the Factorization method with respect to various modelling errors.
\end{abstract}

\section{Introduction}
Electrical Impedance Tomography (EIT) is an imaging technique that allows retrieval of the conductivity distribution inside a body by the application of a current to its 
boundary and measurement of the resulting voltage. The portability and the low cost of electronic devices capable of producing such data makes it an ideal tool for non 
destructive testing or medical imaging. In the last few years,  several imaging devices that use EIT have produced interesting results in the field of 
medical imaging (see \cite{Hol05}) as well as for non destructive testing (see \cite{Yor01}).

The EIT inverse problem has been intensively studied in the mathematical literature since Calder\'on's paper in 1980 (see \cite{Cad80}) 
that formulates the inverse boundary value problem. Then, many authors have been interested in proving uniqueness for the inverse problem and in providing 
efficient algorithms to find the conductivity from boundary measurements (see \cite{MueSil12} and references therein for a complete review on this subject). 
One of these algorithms is the so-called Factorization method, introduced by Kirsch in \cite{Kir98} for locating obstacles from acoustic scattering data and then 
extended by Br\"uhl in \cite{Bru01} to the EIT inverse problem. The main advantage we see in this technique is that it places fewer demands on the data since it 
only locates an embedded inhomogeneity but does not give the conductivity value inside the inclusion. The Factorization method  and its 
regularisation have been studied by Lechleiter \textit{et al.} in \cite{LeHyHa08} in the context of the so-called complete electrode model which was shown by 
Isaacson \textit{et al.} in \cite{ChIsNeGi89} to be close to real-life experimental devices. 

The main purpose of this paper is to show that the Factorization method can be successfully applied to realistic three dimensional EIT inverse problems. To our 
knowledge, the only work presenting three dimensional reconstruction by using the Factorization method in EIT is due to Hanke \textit{et al.} \cite{HanSch08} and 
they focus on the case where the probed domain is the half space. Therefore, we would like to extend this to more complex geometries with inhomogeneous 
background conductivities. To do so, the main difficulty consists in computing the Neumann Green's function of the studied domain since this function is needed to 
apply the Factorization method. In two dimensions with homogeneous background, one can use conformal mapping techniques to obtain the Neumann Green's 
function of many geometries from that of the unit circle (see \cite{BruHan00,HyHaPu07} for example). Clearly, this technique is no longer available in 
 three dimensions and we have to use other ideas. We choose to follow a classical idea presented in \cite{Green}, and more recently used in \cite{HadMig13} in the EIT 
context in two dimensions, because it can be extended to three dimensional problems  and allows to compute Neumann Green's functions for non constant 
conductivities.  It consists in splitting the Green's function into a singular part (which is known) plus a regular part and to compute the regular part which is the solution to 
a well posed boundary value problem with numerical methods such as finite elements or boundary elements. 

We show with various numerical experiments and for different geometries that this technique can be successfully applied in three dimensions to obtain EIT images from 
simulated data.  To make our experiments more realistic, the simulated data we used were produced by using the complete electrode 
model with a limited number of electrodes which covered a part 
of the inspected domain, as it is the case in many applications.
We first compare the results obtained by using the Neumann Green's function of the searched domain and the one of the free 
space. Then we study the influence of various parameters on the quality of the reconstructions and
we conclude our numerical evaluation of the method by studying 
the influence on the reconstructions  of various modelling errors such as errors in electrode positions,
in the shape of the probed domain and in the background conductivity. 

 In section \ref{sec:forward} we present briefly the two main models for the direct EIT problem and introduce some notations 
and definitions. 
In section \ref{sec:inverse} we present the inverse problem and state the main theoretical foundation of this paper. Finally, in section \ref{sec:numerics}, we present 
our numerical implementation of the method and in section \ref{sec:images} we 
give numerical reconstructions for different domains obtained with noisy simulated data.

\section{Forward models in electrical impedance tomography}
\label{sec:forward}
\subsection{The continuum forward model}
Let $\Om \subset \bR^3$ be an open bounded set with Lipschitz boundary $\Ga$. Then, according to the continuum model, the electric potential $u$ inside $\Om$ 
produced by the injection of a current $I$ on the boundary $\Ga$ solves
\begin{equation}
\label{eq:fwd}
\begin{cases}
\div(\sigma\nabla u) =0 \qquad {\rm in} \  \Om, \\
\displaystyle\D{u}{\nu} =I \qquad {\rm on} \  \Ga, \\
\displaystyle \int_\Ga u \,ds=0, 
\end{cases}
\end{equation}
where $\nu$ is the unit normal to $\Ga$ directed outward to $\Om$ and the conductivity $\sigma \in L^\infty(\Om)$ is a real valued function such that there exists 
$c>0$ for which
\[
\sigma(x) \geq c>0 \qquad {\rm for\  almost\  all } \ x \in \Om.
\] 
Let us denote $L^2_\d(\Ga):=\{f \in L^2(\Ga) \ |\ \int_\Ga f ds =0\}$.  It is well known that whenever $I\in L^2_\d(\Ga)$, problem \eqref{eq:fwd} has a unique solution $u 
\in H^1(\Om)$. Then, we can define the so-
called Neumann-to-Dirichlet (NtD) map corresponding to the conductivity $\sigma$ by 
\begin{eqnarray*}
\Ld_\sigma \,: \,& L^2_\d(\Ga) \longrightarrow L^2_\d(\Ga) \\
&I  \longmapsto u|_\Ga
\end{eqnarray*}
where $u$ is the unique solution to \eqref{eq:fwd}.

\subsection{The complete electrode model}
A more accurate model, the so-called complete electrode model, takes into account the fact that in practise the current is applied on $\Ga$ through a finite number 
$M \in \mathbb{N}$ of electrodes. Let us denote by $(E_j)_{j=1,\ldots,M}$ these electrodes, for each $j =1,\ldots,M$, $E_j$ is a connected open subset of $\Ga$ 
with Lipschitz boundary and non zero measure. According to the experimental setups, we assume in addition that the distance between two electrodes is strictly 
positive, that is $\overline{E_i} \cap \overline{E_j} = \emptyset$ when $i\neq j$. We also introduce $\mathcal{G}^M:= \Ga \setminus \cup_{i=1}^M E_i$ the gap 
between the electrodes. Following the presentation in \cite{Hyv04}, let us introduce the subspace of $L^2_\d(\Om)$ of functions that are 
constant on each electrode and that vanish on the gaps between the electrodes
\[
T^M_\d:=\left \{f \in L^2_\d(\Ga) \ | \ f = \sum_{i=1}^M \chi_{E_i} f_i,\ f_i \in \bR,\ 1\leq i \leq M \right\}
\]
where $\chi_V$ stands for the characteristic function of some domain $V$. For simplicity, in the following we will not make the distinction between 
an element $f$ of $T^M_\d$ and the associated vector $(f_i)_{i=1,\ldots, M}$ of $\bR^M$. Let us denote by $I \in T^M_\d$ the injected current which is constant and 
equal to $I_i$ on electrode $i$. The voltage potential $u$  then solves
\begin{equation}
\label{eq:CEM}
\begin{cases}
\div(\sigma\nabla u) =0 \qquad {\rm in} \  \Om, \\
\frac{1}{|E_i|}\displaystyle\int_{E_i}\sigma\D{u}{\nu} ds=I_i \qquad {\rm for} \  i=1,\ldots M, \\
u+z \sigma \D{u}{\nu} =U_i\qquad {\rm for}\ x \in E_i, \  i=1,\ldots M,\\
\sigma \D{u}{\nu} =0  \qquad {\rm for}\ x \in \Ga \setminus\cup_{i=1}^M \overline{E_i}
\end{cases}
\end{equation}
where $z \in L^2(\Ga)$ is the so-called contact impedance and $U\in T^M_\d$ is the unknown measured voltage potential 
on the electrodes. We assume that the contact impedance satisfies $z(x)\geq c>0$ for almost all $x\in \Ga$.

For all $I \in T^M_\d$ there exists a unique $(u,U) \in H^1(\Om)\times T^M_\d$ that solves problem \eqref{eq:CEM} (see \cite{SoChIs92} for more details about the model 
and the proof of existence and uniqueness) and this solution depends 
continuously on $I$. Similarly to the continuum model, we define the (finite dimensional) Neumann-to-Dirichlet map associated with the conductivity $\sigma$ 
by 
\begin{eqnarray*}
\Sigma^M_\sigma \,: \,& T^M_\d \longrightarrow  T^M_\d \\
&I  \longmapsto U
\end{eqnarray*}
where $(u,U)$ is the unique solution to \eqref{eq:CEM}.

\section{The inverse problem}
\label{sec:inverse}
\subsection{Statement of the inverse problem}
 Let us denote by $\sigma\in L^\infty(\Om)$ the conductivity in the presence of an inclusion. We assume that $\sigma(x)\geq c>0$ for almost all $x\in \Om$ and that 
there exists a domain $D\subset \Om$  such that
\[
\sigma(x) = \sigma_0(x) + \ga(x) \chi_D(x)
\]
where $\ga \in L^\infty(D)$ is either a  positive or a negative function on $D$ and  $\sigma_0 \in C^{0,1}(\Om)$ is the known conductivity of the background and is such that $\sigma_0(x) \geq c>0$ for  almost  all  $ x \in \Om$. We assume moreover that $\overline{D}\subset \Om$ and that $\Om\setminus
\overline{D}$ is connected with Lipschitz boundary. The inverse problem we treat consists in finding the indicator function $\chi_D$ from the knowledge of the finite 
dimensional maps $\Sigma_{\sigma_0}^M$ and $\Sigma_\sigma^M$. 

As stated in the next section (see Theorem \ref{th:factocont}), the Factorization method
 provides an explicit  formula of the indicator function $\chi_D$ from the knowledge of $\Lambda_{\sigma_0}$ and $\Lambda_{\sigma}$. Nevertheless, as 
we explain later on,  $\Sigma_{\sigma_0}^M$ and $\Sigma_{\sigma}^M$  actually give a good approximation of the characteristic function of $D$ provided that the 
number of electrodes $M$ is sufficiently large and that they cover most of $\Ga$.

\subsection{The Factorization method to find the support of an inclusion}
 For each point $z$ in $\Omega$, let us define the Green's function for the background problem 
with Neumann boundary conditions $N(\cdot,z) \in L_\d^2(\Om)$ which is the solution to 
\begin{equation}
\label{eq:Green}
\begin{cases}
\div_x(\sigma_0 (x) \nabla_x N(x,z)) =-\delta_z \qquad \rm{in}\, \Om, \\
\displaystyle \sigma_0 \D{N(x,z)}{\nu(x)} = -\frac{1}{|\Ga|} \qquad \rm{on} \, \Ga, \\
\displaystyle \int_\Ga N(x,z) ds(x) =0,
\end{cases}
\end{equation}
where $\delta_z$ is the Dirac distribution at point $z\in \Om$. Then,
\begin{equation}
\label{eq:phizd}
\phi_z^d(x) := d \cdot \nabla_z N(x,z)
\end{equation}
is the electric potential created by a dipole located at point $z \in \Om$ of direction $d\in \bR^3$ such that $|d|=1$. 

We use these dipole test functions in the next Theorem (which is a slightly reformulated version of  \cite[Proposition 4.4]{Bru01})  in the expression of the 
characteristic function of the inclusion $D$.
\begin{theorem}
\label{th:factocont}
Take $d\in \bR^3$ such that $|d|=1$ and let us denote $(\ld_i,\psi_i)_{i=1,\ldots ,+\infty}$ the eigenvalues and eigenvectors of the self adjoint and compact 
operator $\Ld_{\sigma_0}-\Ld_\sigma$. Then the characteristic function $\chi_D$ of $D$ is given by
\begin{equation}
\label{eq:caracPicard}
\chi_D(z) =
{\rm sgn}  \left[\left(\sum_{i=1}^{+\infty} \frac{(\phi_z^d,\psi_i)^2_{L^2(\Ga)}}{|\ld_i|}\right)^{-1}\right].
\end{equation}
\end{theorem}
\begin{remark}
\label{re:equiv}
An equivalent (and maybe more classical) statement of Theorem \ref{th:factocont} is the following range test
\[
z \in D \quad \Longleftrightarrow \quad \exists \, g_z \in L^2_\d(\Ga) \textrm{ such that } |\Ld_{\sigma_0}-\Lambda_\sigma|^{1/2}g_z = \phi^d_z
\]
for $d$ a given unit vector of $\bR^3$. Then, \eqref{eq:caracPicard} is given by the inverse $L^2$ squared norm of $g_z$  that solves $|\Ld_{\sigma_0}-\Lambda_
\sigma|^{1/2}g_z = \phi^d_z$  and which is finite if and only if this equation has a solution in $L^2_\d(\Ga)$.
\end{remark}
In \cite{LeHyHa08} the authors state a convergence result that justifies the use of the finite dimensional  map $\Sigma_{\sigma_0}^M-\Sigma_{\sigma}^M$  instead 
of  $\Ld_{\sigma_0}-\Ld_\sigma$ to obtain an approximation of  the characteristic function of $D$. We recall here briefly the main result 
they obtain. Let us 
introduce $P^M : L^2_\d(\Ga) \rightarrow T^M_\d$ the $L^2(\Ga)$ orthogonal projector defined by
\[
P^M f = \sum_{i=1}^M\chi_{E_i}\frac{1}{|E_i|}\left(\int_{E_i} f ds+ \frac{1}{M}\int_{\mathcal{G}^M} f ds\right)
\]
and let us denote by $\Sigma^M :=(\Sigma_{\sigma_0}^M-\Sigma_{\sigma}^M) : T^M_\d \rightarrow T^M_\d$ the difference data operator. Then, if the number of 
electrodes $M$ goes to infinity 
and if the gap between the electrodes decreases sufficiently fast, we deduce from \cite[Theorem 8.2]{LeHyHa08} that for every $M$ there exists a truncation index 
$0<R(M)< M$ such that for a given $z \in \Om$ the sequence
\begin{equation*}
\frac{1}{\|P^M\phi_z^d\|^2_{L^2(\Ga)}}\sum_{i=1}^{R(M)} \frac{(P^M \phi_z^d,\psi^M_i)^2_{L^2(\Ga)}}{|\ld^M_i|}
\end{equation*}
converges when $M$ goes to infinity if and only if $z$ is in $D$. In this expression, $(\ld^M_i,\psi^M_i)_{i=1,\ldots ,M}$ are the eigenvalues and eigenfunctions of 
the finite dimensional map $\Sigma^M$. In practice, it is not easy to determine the truncation index $R(M)$ from the data (see \cite{BruHan03} for a heuristic 
method) but since we will consider the case where 
we have few measurements, we will see that no regularisation is actually needed and we will take $R(M)=M-1$. This result tells us that 
the function 
\[
\chi^M_d(z):=\|P^M\phi_z^d\|^2_{L^2(\Ga)}\left(\sum_{i=1}^{R(M)} \frac{(P^M\phi_z^d,\psi^M_i)^2_{L^2(\Ga)}}{|\ld^M_i|}\right)^{-1}
\]
should be an approximation of the indicator function of $D$ in the sense that $\chi^M_d$ should be greater inside $D$ than outside $D$.
\begin{remark}
\label{re:rangetest}
As in the continuous setting (see Remark \ref{re:equiv}), for each point $z$ the characteristic function can be defined by using the inverse of the squared $L^2$ norm of the solution $g_z^M$ to the linear equation
\[
| \Sigma^M|^{1/2} g^M_z = P^M \phi^d_z.
\]
\end{remark}

\section{Numerical implementation}
In this section we give a precise definition of the data we used for inversion and we present our implementation of the Factorization method
and of the computation of the dipole test functions. 
\label{sec:numerics}
\subsection{Data sets and numerical implementation of the indicator function}
In practise, one does not know the full and noiseless map $\Sigma^M$ but only the map $\hat \Sigma^M : V \rightarrow T_\d^M$ where $V \subset T_\d^M$ is 
the 
set of currents that one injects in the tested object and for all $I\in V$ we have 
\[
\hat \Sigma^M I = \Sigma^M I + \epsilon
\]
where $\epsilon$ denotes the noise in the data.  As a consequence, the data consist of the $M\times N$ matrix $\hat \Sigma^M$ 
where $N$ is the dimension of $V$. Each column of this matrix is a vector that contains the voltage at each electrode.  In the following, we will consider the case of synthetic data.  
To produce them, we compute a noiseless map $\wt \Sigma^M$  by using finite elements implemented with the software FreeFem++ (see \cite{ff++}) to solve 
equations \eqref{eq:CEM} for each current $I$ in the admissible set of currents $V$. We refer to \cite{VaVaSaKa99} for more details on the variational formulation we 
use to solve \eqref{eq:CEM}. Then, we build $\hat \Sigma^M$ by adding artificial noise to $\wt \Sigma^M$ that is
\[
\hat \Sigma^M = \wt \Sigma^M  + \eta \mathcal{N}\cdot \wt\Sigma^M
\]
where $ \mathcal{N}$ is a matrix of size $M\times N$ whose $(i,j)$ element is a random number generated with a normal distribution  and $\eta$ is a real number 
chosen such that 
\begin{equation}
\label{eq:delta}
\delta := \frac{\| \wt \Sigma^M-\hat \Sigma^M\|}{\|\wt\Sigma^M\|}
\end{equation}
is a given level of noise. The $\cdot$ denotes the term by term multiplication between two matrices  and $\|\cdot\|$ denotes the Frobenius norm.

Let us introduce the singular values $(\hat \sigma^M_i)_{i=1,\ldots N}$ and singular vectors $(\hat u^M_i,\hat v^M_i)_{i=1,\ldots,N}$ of  $\hat \Sigma^M$ that satisfy
\[
\hat \Sigma^M \hat v^M_i= \hat \sigma^M_i \hat u^M_i \quad {\rm and} \quad (\hat \Sigma^M)^T \hat u^M_i= \hat \sigma^M_i \hat v^M_i 
\]
for all $i=1,\ldots,N$. Then, the function $\chi^M_d$ introduced in the previous section corresponds to the inverse of 
\[
f^M(z,d):=\frac{1}{\sum_{i=1}^{R(M)}(P^M\phi_z^d,\hat u^M_i)^2}\sum_{i=1}^{R(M)} \frac{(P^M\phi_z^d,\hat u^M_i)^2}{\hat \sigma^M_i}
\]
 since 
\[
g^M_z:=\sum_{i=1}^{R(M)} \frac{(P^M\phi_z^d,\hat u^M_i)}{(\hat \sigma^M_i)^{1/2}} \hat v_i^M
\]
 solves $|\hat \Sigma^M|^{1/2} g^M_z = P^M \phi^d_z$  (see Remark \ref{re:rangetest}). As mentioned before,  we take the truncation index $R(M)$ as large as we can which is equal to the 
dimension of the range of $\Sigma^M$. Finally, to limit artefacts, we will use the function
\[
{\rm Ind}(z) :=\left(\sum_{d \in S} f^M(z,d) \right)^{-1}
\]
as an indicator function of $D$ where $S$ is a set of $8$ unit vectors of $\bR^3$. The choice of the number of dipole directions is based on  experimental observations.

\subsection{Computation of the dipole potentials}
In three dimensions, the Green's function of the Laplace equation in the free-space is given by
\[
\Phi_z(x) = \frac{1}{4\pi}\frac{1}{|x-z|^3}.
\]
and for any unit vector $d\in \bR^3$ and points $x\neq z$ we define the associated dipole potential
\[
\hat \phi_z^d(x) :=d\cdot \nabla_z \Phi_z(x) =  -\frac{1}{2\pi}\frac{d\cdot(x-z)}{|x-z|^3}.
\]
We remark that the image principle used in two dimensions  to compute the Neumann Green's function for a circle is not valid
 anymore in three dimensions for a sphere. 
 Nevertheless (see \cite{Green,HadMig13} for example), for   $\sigma_0 \in C^{0,1}(\Om)$ and  for $z, x \in \Om$ we can decompose $\phi_z^d$ as 
\[
\phi_z^d(x) = V_z^d(x) + \hat \phi_z^d(x,\sigma_0(z))
\]  
where
\[
\hat \phi_z^d(x,\sigma_0(z)) := d\cdot \nabla_z \left( \frac{\Phi_z(x) }{\sigma_0(z)} \right)
\] 
and $V_z^d(x)$  solves the following conductivity problem:
\begin{equation}
\label{eq:GreenRegul}
\begin{cases}
\div_x(\sigma_0 (x) \nabla_x V_z^d) =\rm{div}\left[(\sigma_0(z) - \sigma_0(x)) \nabla_x \hat \phi_z^d(x,\sigma_0(z))\right] \qquad \rm{in}\, \Om, \\
\displaystyle \sigma_0(x) \D{V_z^d(x)}{\nu(x)} = - \sigma_0(x) \D{\hat \phi_z^d(x,\sigma_0(z))}{\nu(x)}\qquad \rm{on} \, \Ga, \\
\displaystyle \int_\Ga V_z^d(x) ds(x) =\int_\Ga  \hat \phi_z^d(x,\sigma_0(z))ds(x).
\end{cases}
\end{equation}
One can actually compute $\nabla_x \hat \phi_z^d(x,\sigma_0(z))$ and realise that this function is singular for $x=z$ but in $L^2(\Om)$. Therefore, by elliptic regularity we deduce that $V_z^d(x)$ is in 
$H^1(\Om)$ and  we compute it by using the finite element software FreeFem++ to obtain a good approximation of $\phi_z^d(x_i)$ for each electrode 
position $x_i$ and each point $z$ in $\Om$. In the case of an homogeneous background, 
we have verified the accuracy of our approximation by comparing the 
finite element solution  $V_z^d(x)$ of \eqref{eq:GreenRegul} with a boundary element solution 
computed with the software BEM++ (see \cite{bempp}). In the case of a homogeneous background ($\sigma_0$ is constant) there is no source
term inside $\Omega$ in equations \eqref{eq:GreenRegul}. Therefore, as it has been observed   in \cite{BruHan00}, one can first compute 
the Neuman-to-Dirichlet map  $\Lambda_{\sigma_0}$ associated with the continuum model and obtain the solutions to
equation \eqref{eq:GreenRegul} with a simple matrix vector product.

In the following, we will study the impact  of using the zero average dipole of the free space: $ \wt \phi_z^d(x):= 
\hat \phi_z^d(x)  - 1/|\Ga|\int_{\Ga} \hat \phi_z^d(x)$ on the quality of the reconstructions.
To this end, we introduce the new indicator function
\[
\wt {\rm Ind}(z) :=\left(\sum_{d \in S} \wt f^M(z,d) \right)^{-1}
\]
where $\wt f^M $ corresponds to $f^M $ with $\phi_z^d$ replaced by $\wt \phi_z^d$.

\section{Numerical simulations and error analysis}
\label{sec:images}
In this section, we show that the Factorization method successfully applies to three dimensional imaging problems
 in rather complicated geometries and with partial 
covering of the boundary $\Ga$ with electrodes. We also test the sensibility of the method with respect to various 
experimental errors such as errors in the shape of the domain, errors in the background conductivity  
and errors in the electrode's placement. 
In what follows we choose to  not show plots of the indicator functions but plots of an 
iso-surface of the indicator functions $\rm \wt {Ind}$ and $\rm Ind$ 
in red on the images.  The choice of  the 
iso-surface is a complicated question and to our knowledge there is no systematic way to do it; we arbitrarily
 choose to show the iso-surface of value $0.9$ for the indicator function which we normalise between $0$ and $1$.  
 We show in Figure \ref{fig:headtruncation} how this parameter influences the quality of the reconstruction.
To overcome this known difficulty, a solution would be to use the indicator function obtained with the Factorization as an initial 
guess for a level set approach 
as it is done in \cite{Nic13} in the context of acoustic scattering or use it for regularisation
of a linear inversion method as it is proposed in \cite{ChHaSepre}.

\subsection{Preliminary experiments}
\label{sec:prelim}
\begin{figure}
\centering
\subfigure[Domain  $\Om$ (in blue) and electrodes (coloured squares).]{\raisebox{.5cm}{\includegraphics[width=.22\textwidth]{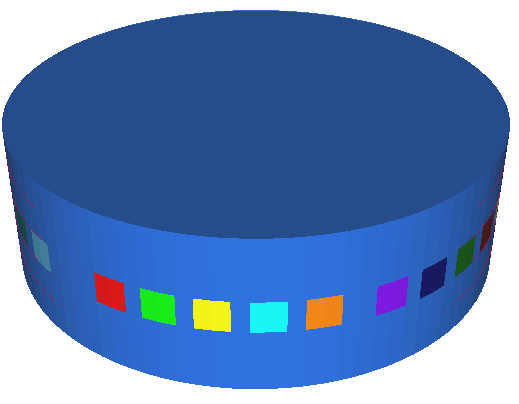}}\label{fig:cyl2setup}}\hfill
\subfigure[One inclusion using $\phi_z^d$; $E_c = 0.12$.]{\includegraphics[width=.34\textwidth]{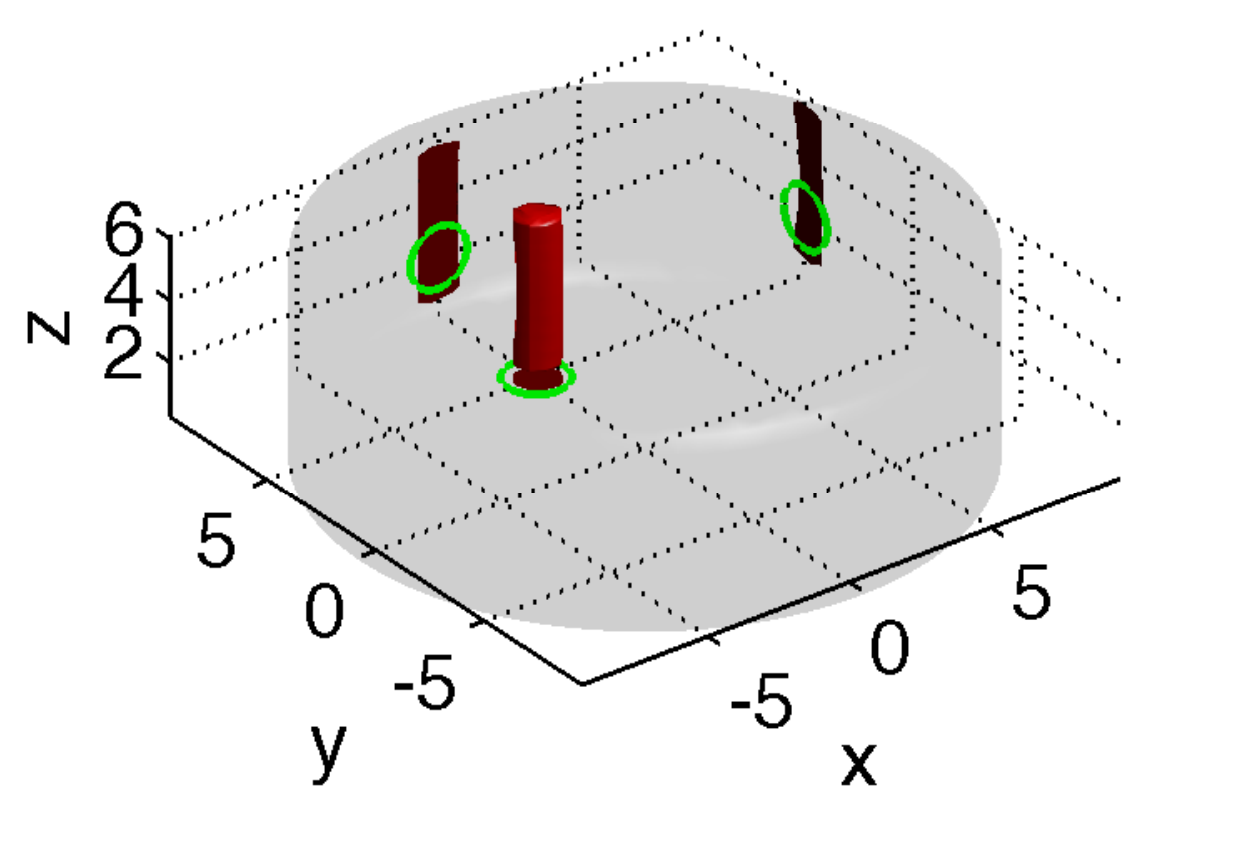}\label{fig:cyl21}} \hfill
\subfigure[One inclusion using $\tilde\phi_z^d$; $\wt E_c = 0.18$.]{\includegraphics[width=.34\textwidth]{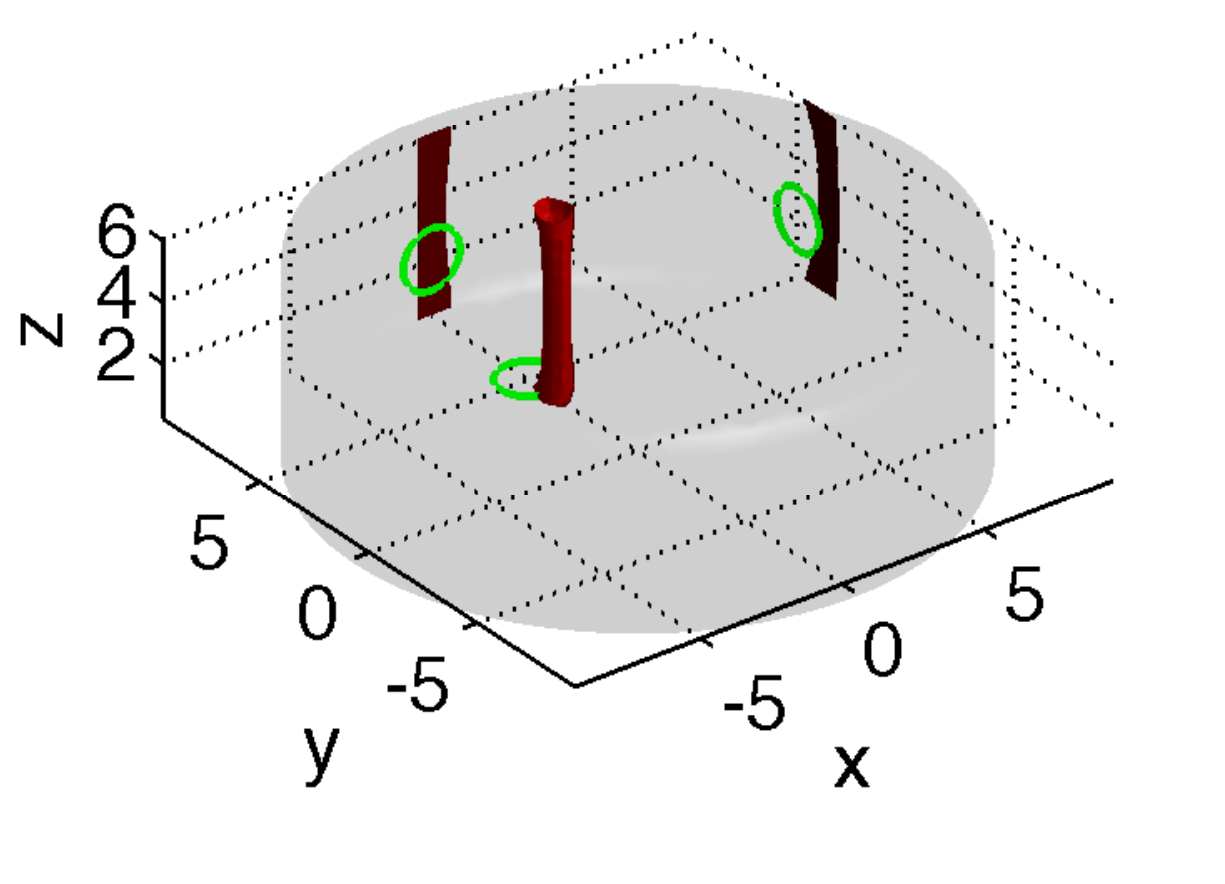}\label{fig:cyl21b}}\hfill
\subfigure[Two inclusions using $\phi_z^d$; $E_c = 0.13$.]{\includegraphics[width=.34\textwidth]{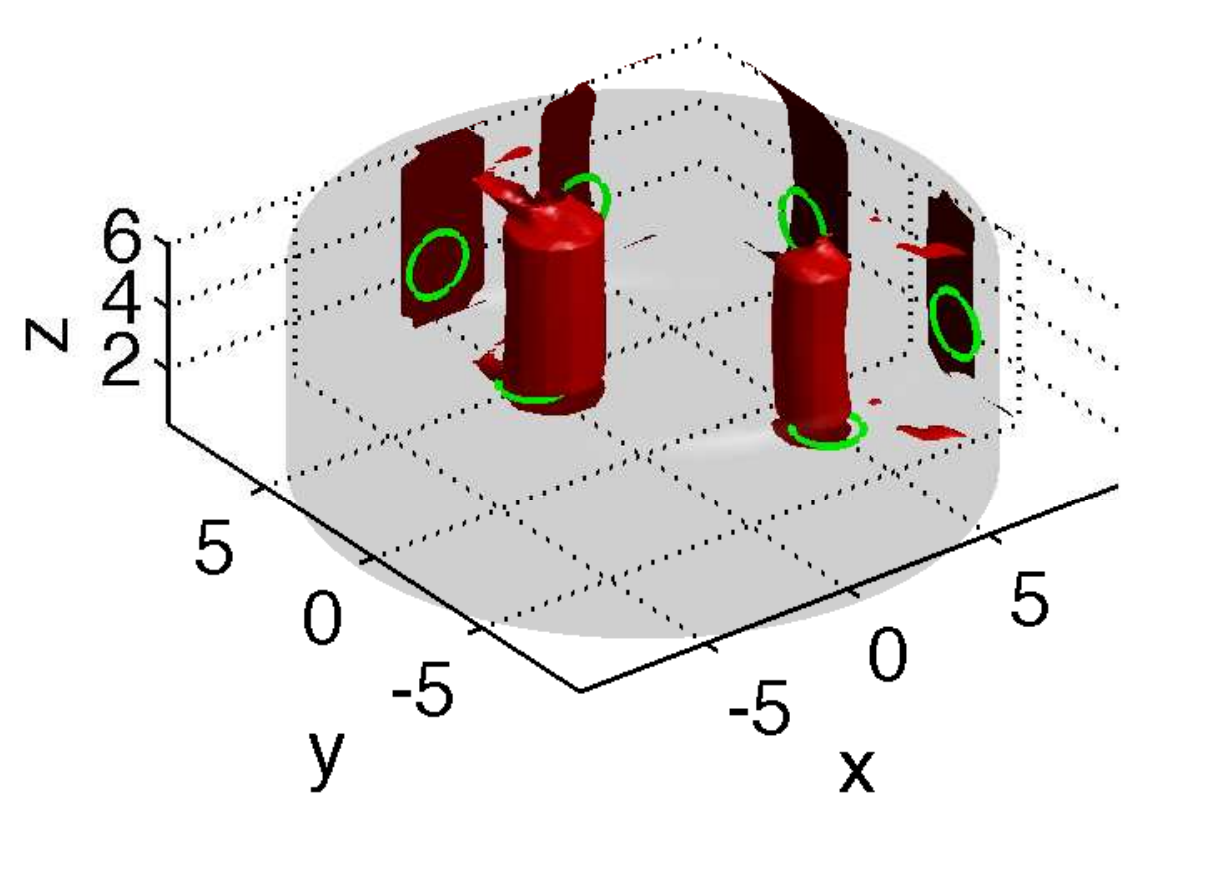}\label{fig:cyl22}}\hspace{1cm}
\subfigure[Two inclusions using $\wt \phi_z^d$; $\wt E_c = 0.21$.]{\includegraphics[width=.34\textwidth]{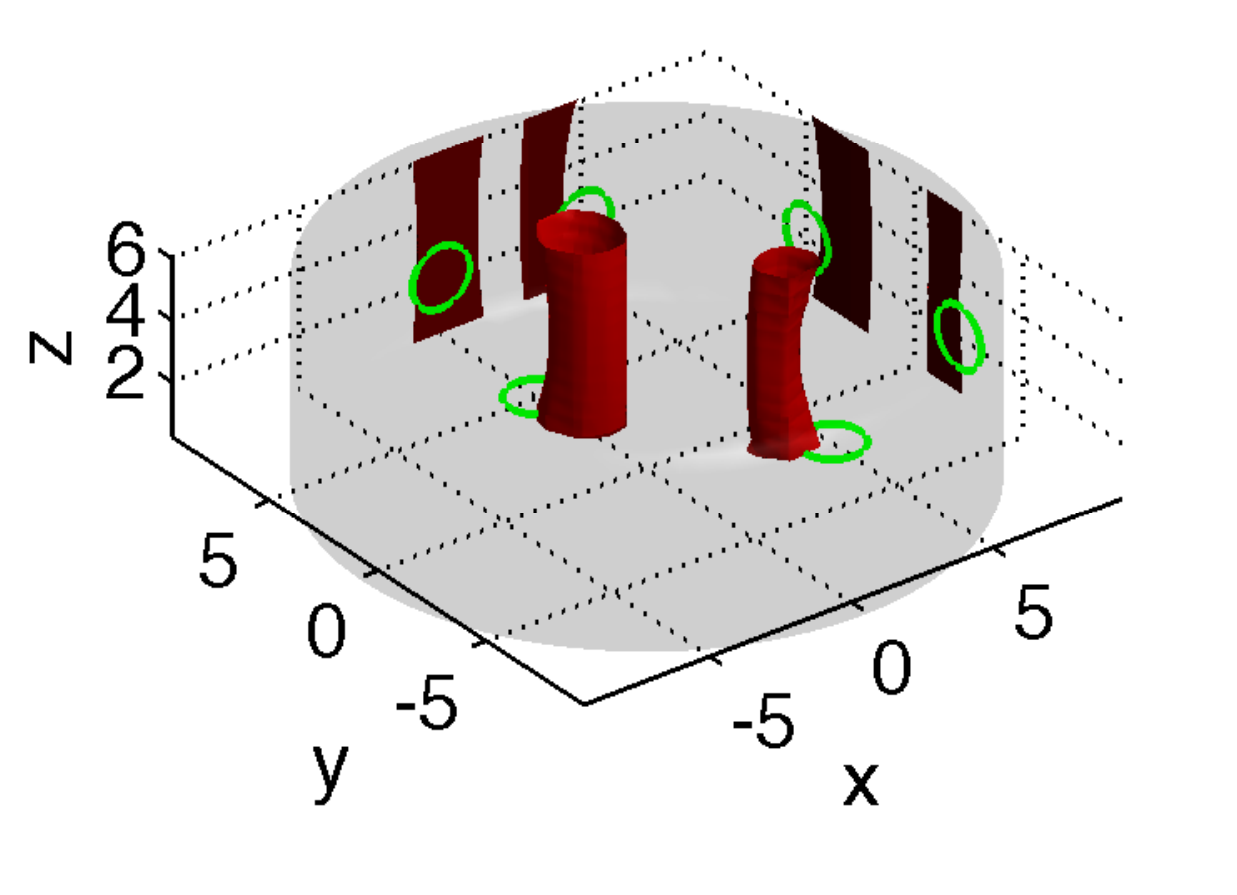}\label{fig:cyl22b}}
\caption{Reconstruction for a cylindrical geometry with one ring of electrodes by adding $\delta=1\%$ of noise to synthetic the data. The geometry of the domain is shown in Figure (a) and in Figures (b)-(e) we represent reconstructions for different cases. In red we plot the iso surface of value $0.9$ of the indicator function as well as its projection on the border of the plotting domain and in green is the projection of the true object. }
\label{fig:cylinder1}
\end{figure}

\begin{figure}[h!]
\centering
\subfigure[Domain  $\Om$ (in blue) and electrodes (coloured squares).]{\raisebox{.5cm}{\includegraphics[width=.22\textwidth]{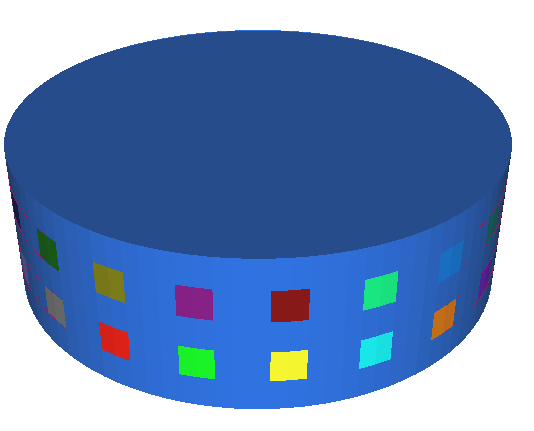}\label{fig:cyl3setup}}}\hfill \hspace{.4cm}
\subfigure[One inclusion using $\phi_z^d$; $E_c = 0.03$.]{\includegraphics[width=.34\textwidth]{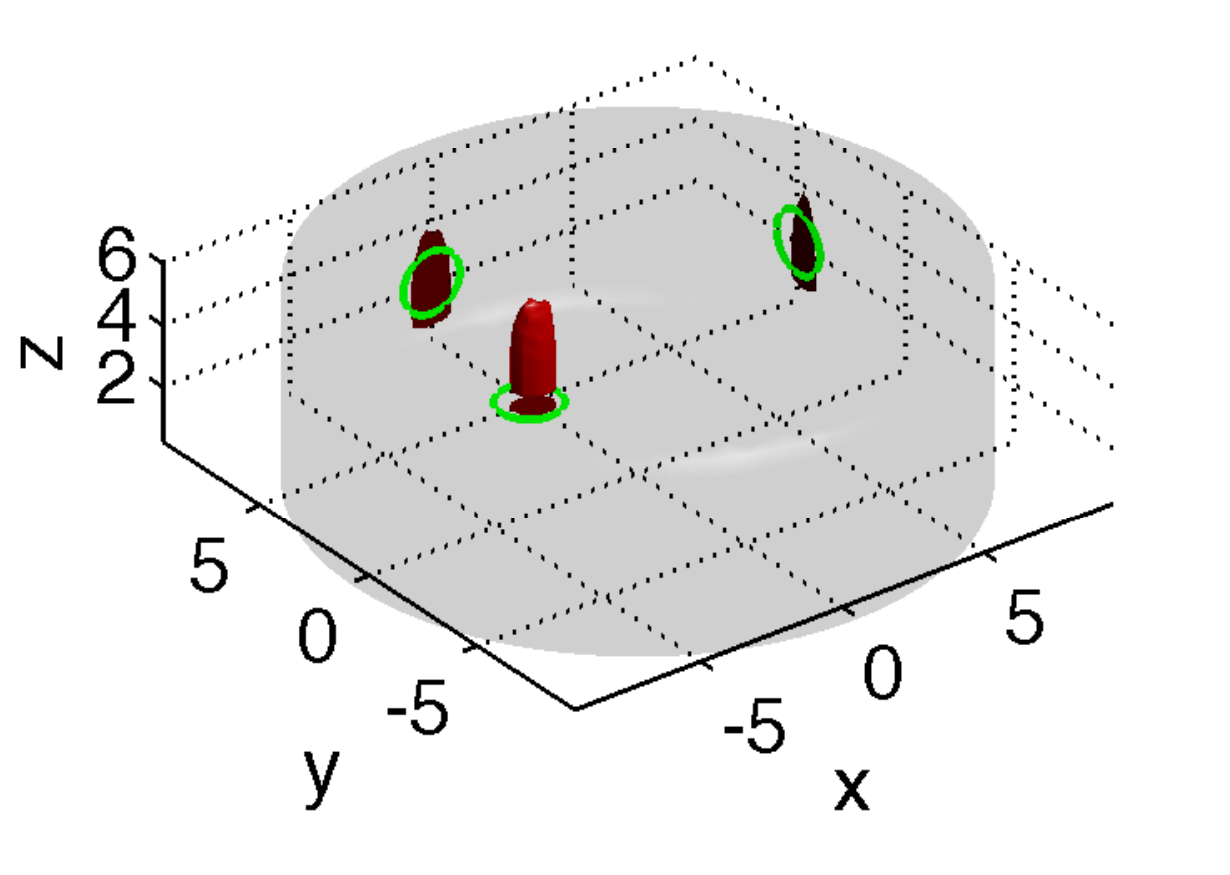}\label{fig:cyl31}} \hfill
\subfigure[Two inclusions using $\phi_z^d$; $E_c = 0.06$.]{\includegraphics[width=.34\textwidth]{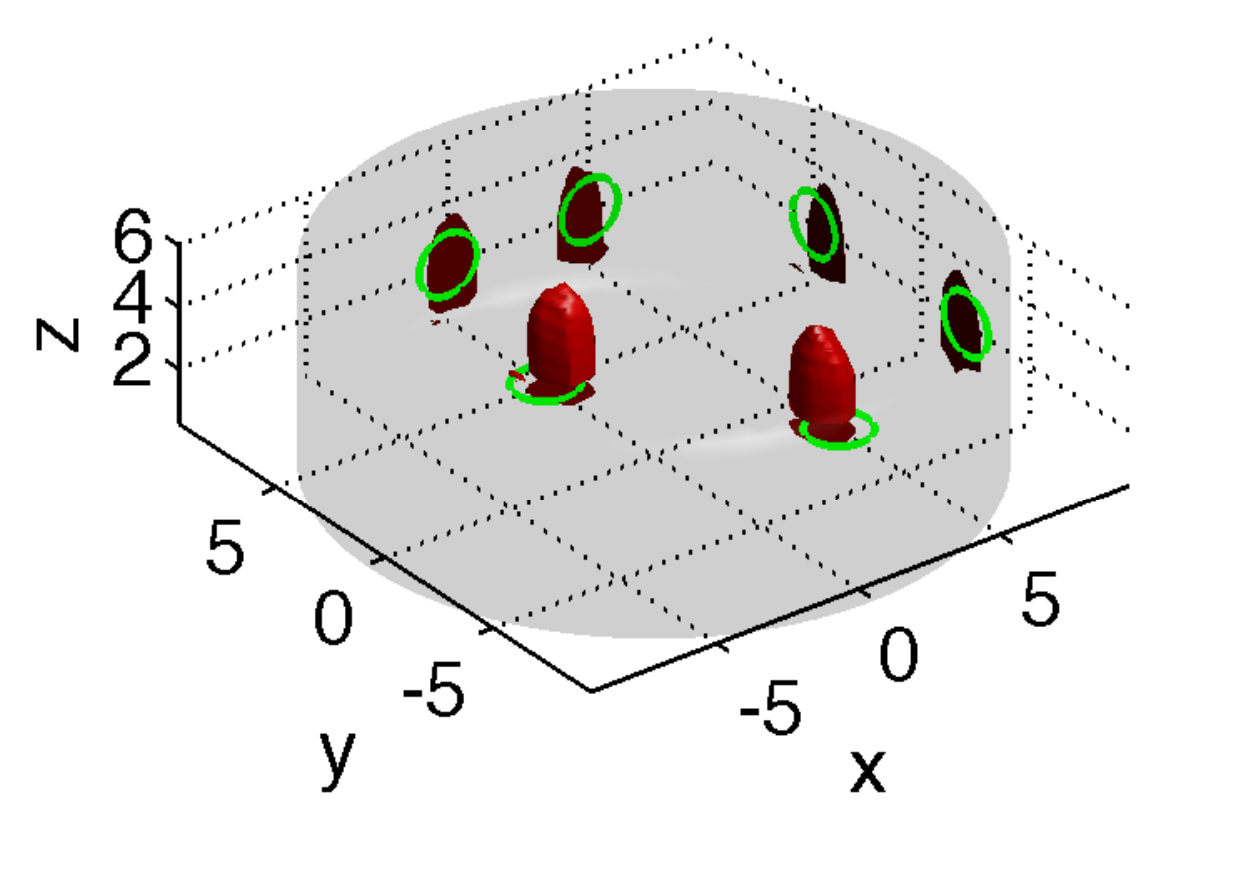}\label{fig:cyl32}}\hspace{1cm}
\caption{Reconstruction for a cylindrical geometry with two rings of electrodes by adding $\delta=1\%$ of noise to the synthetic data (same colours as in Figure \ref{fig:cylinder1}). }
\label{fig:cylinder2}
\end{figure}

\begin{figure}
\centering
\subfigure[Domain  $\Om$ (in blue) and electrodes (coloured circles).]{\includegraphics[width=.22\textwidth]{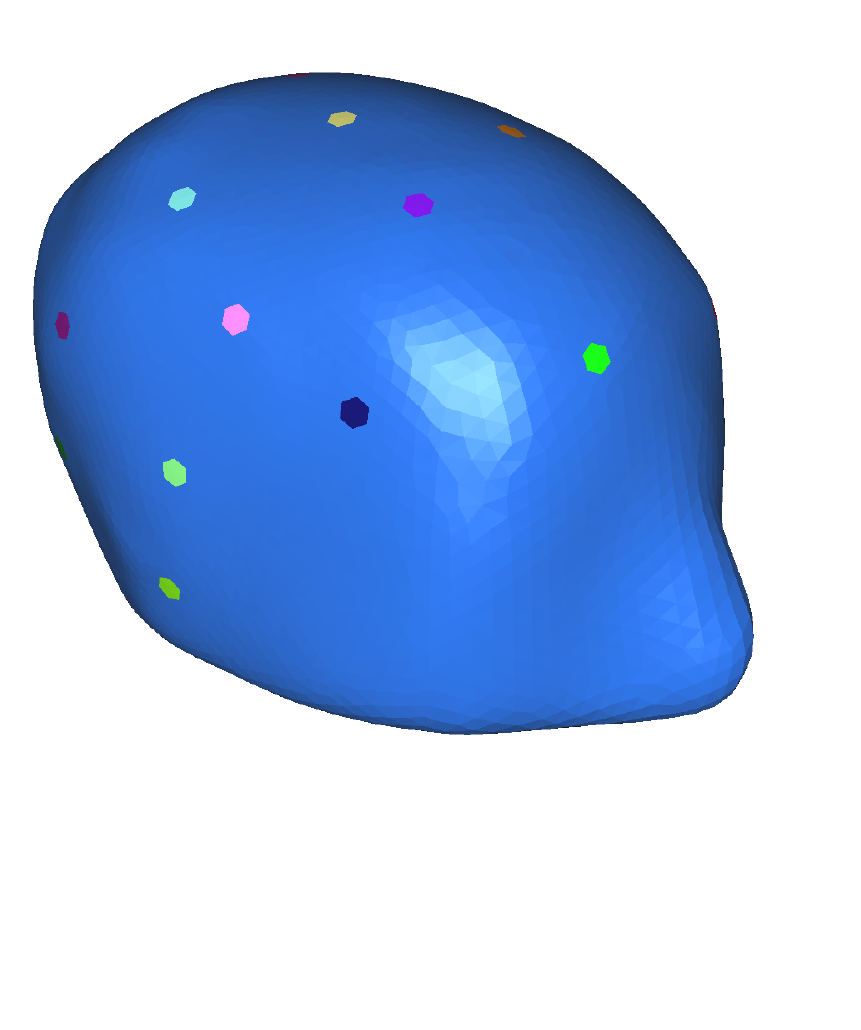}\label{fig:headhmosetup}}\hfill
\subfigure[ One inclusion in the middle using $\phi_z^d$; $E_c = 0.007$.]{\includegraphics[width=.34\textwidth]{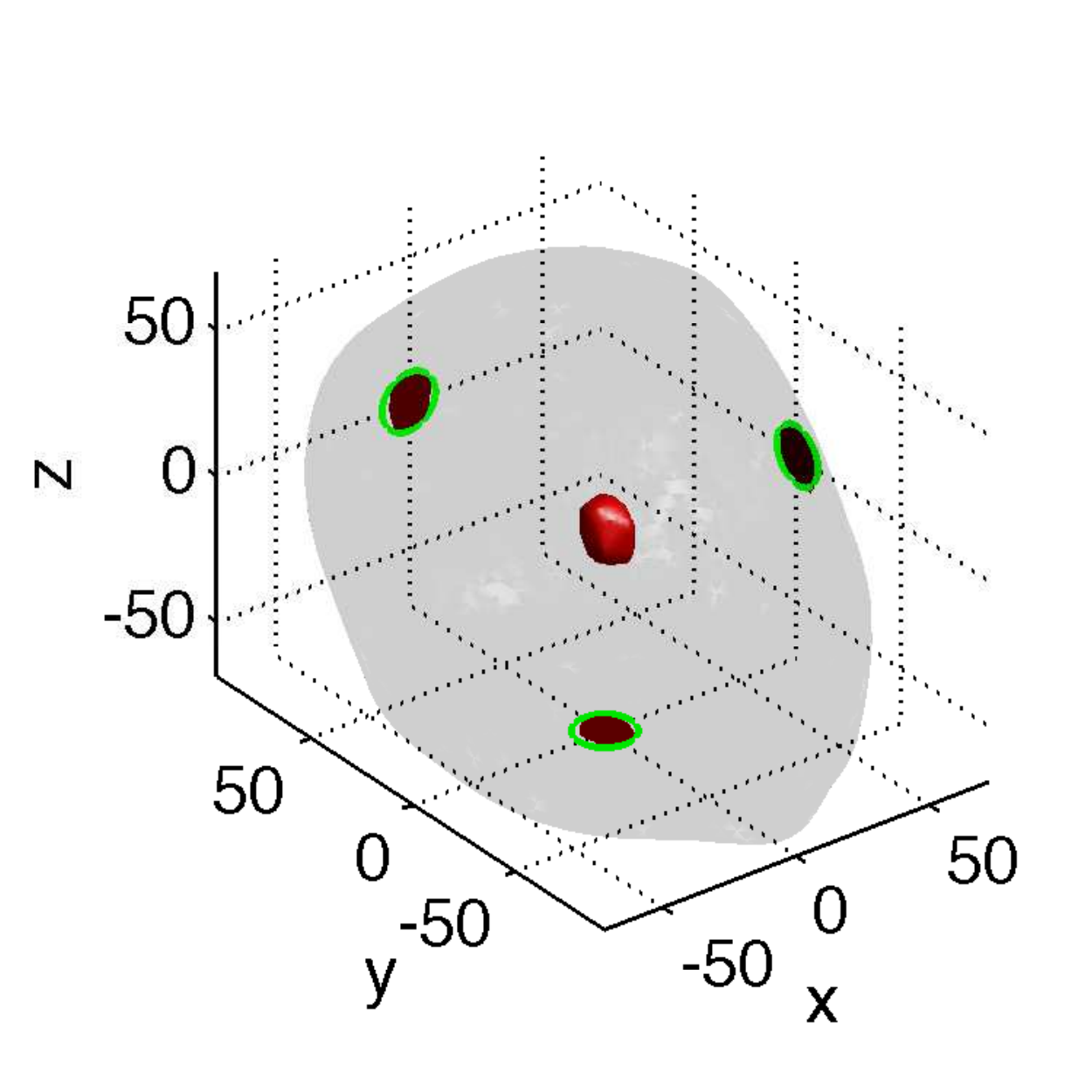}\label{fig:hmiddle}} \hfill
\subfigure[One inclusion in the  back side using $\phi_z^d$; $E_c = 0.01$.]{\includegraphics[width=.34\textwidth]{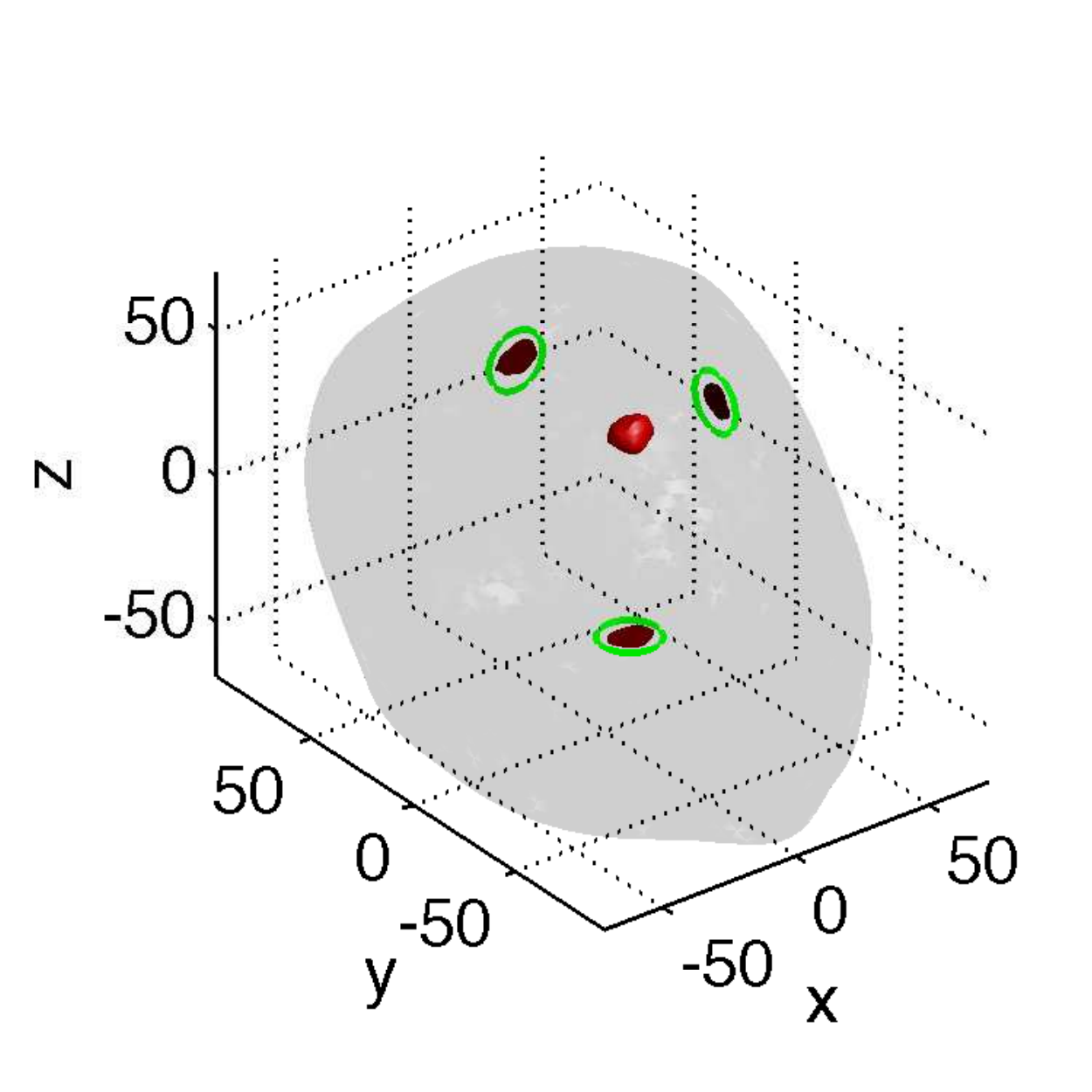}\label{fig:hback}} 
\caption{Reconstruction for a head shape geometry with homogeneous background by adding $\delta=1\%$ of noise to the data (same colours as in Figure \ref{fig:cylinder1}).}\label{fig:headhomogeneous}
\end{figure}

\begin{figure}
\centering
\subfigure[Inhomogeneous layered structure.]{\includegraphics[width=.22\textwidth]{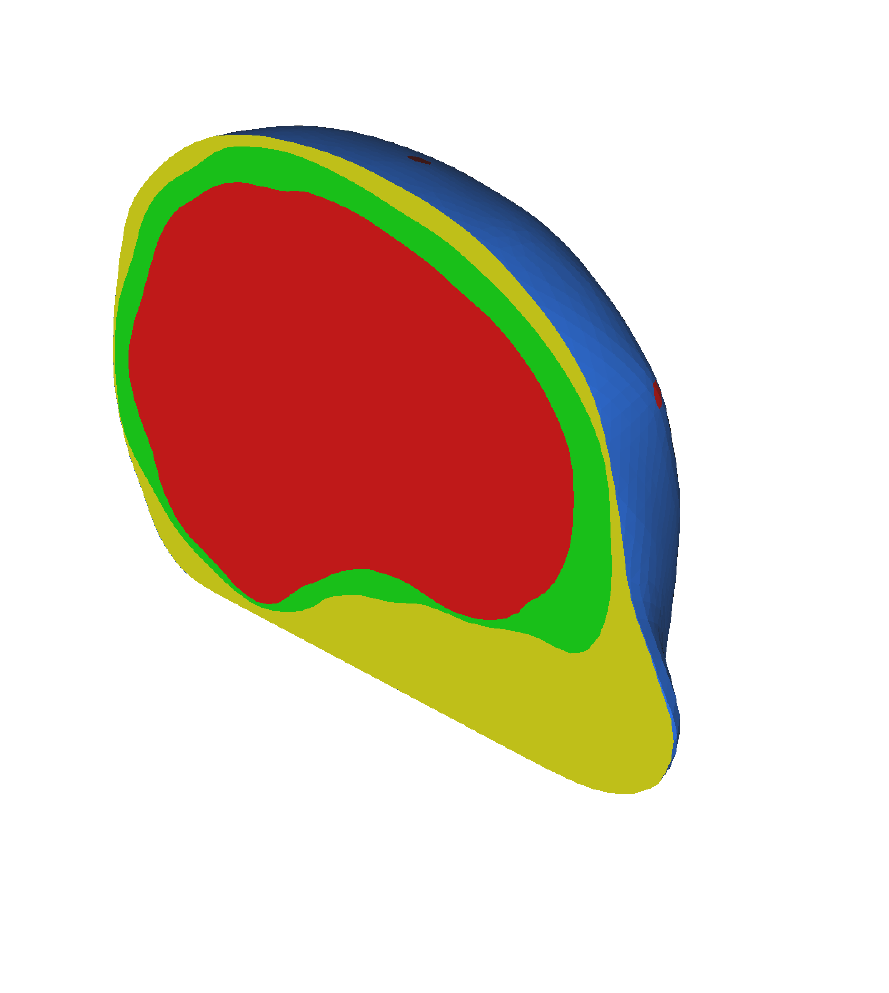}\label{fig:skullsetting}}\hfill
\subfigure[One inclusion in the middle using $\phi_z^d$; $E_c = 0.01$.]{\includegraphics[width=.34\textwidth]{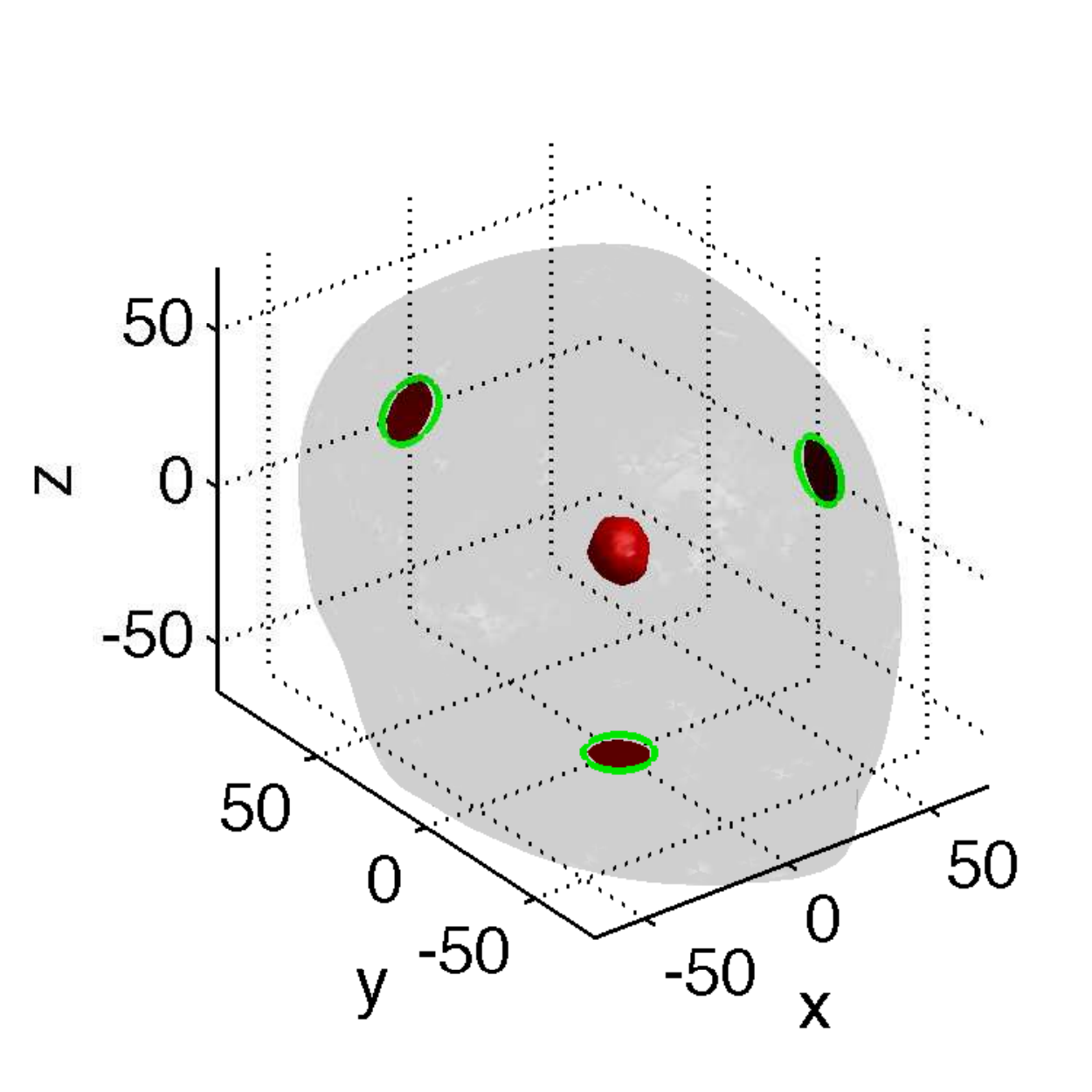}\label{fig:hsmiddle}}\hfill
\subfigure[One inclusion in the middle using $\wt \phi_z^d$; $\wt E_c = 0.08$.]{\includegraphics[width=.34\textwidth]{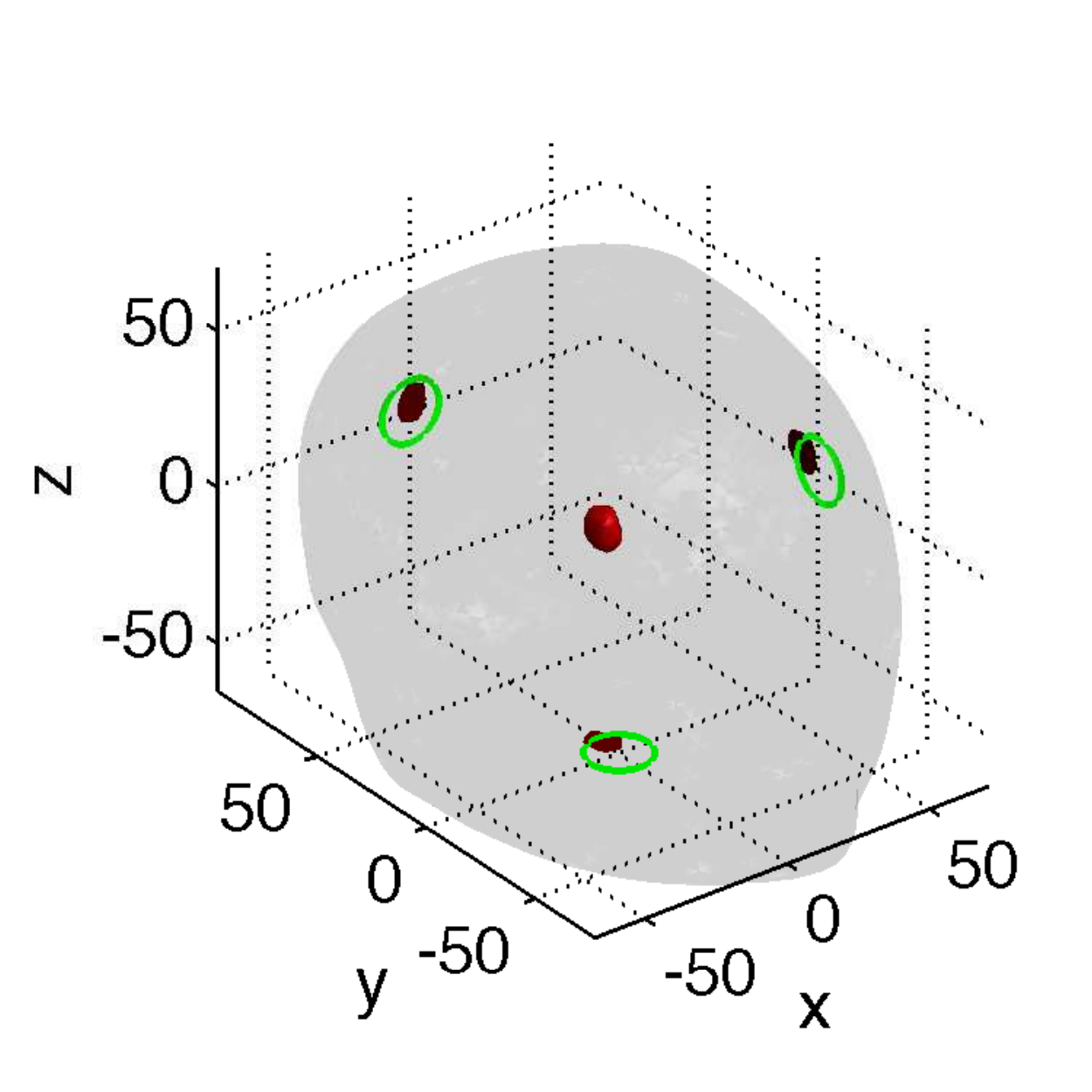}\label{fig:hsmiddleb}}\hfill
\subfigure[One inclusion in the bask side using $\phi_z^d$; $E_c = 0.01$.]{\includegraphics[width=.34\textwidth]{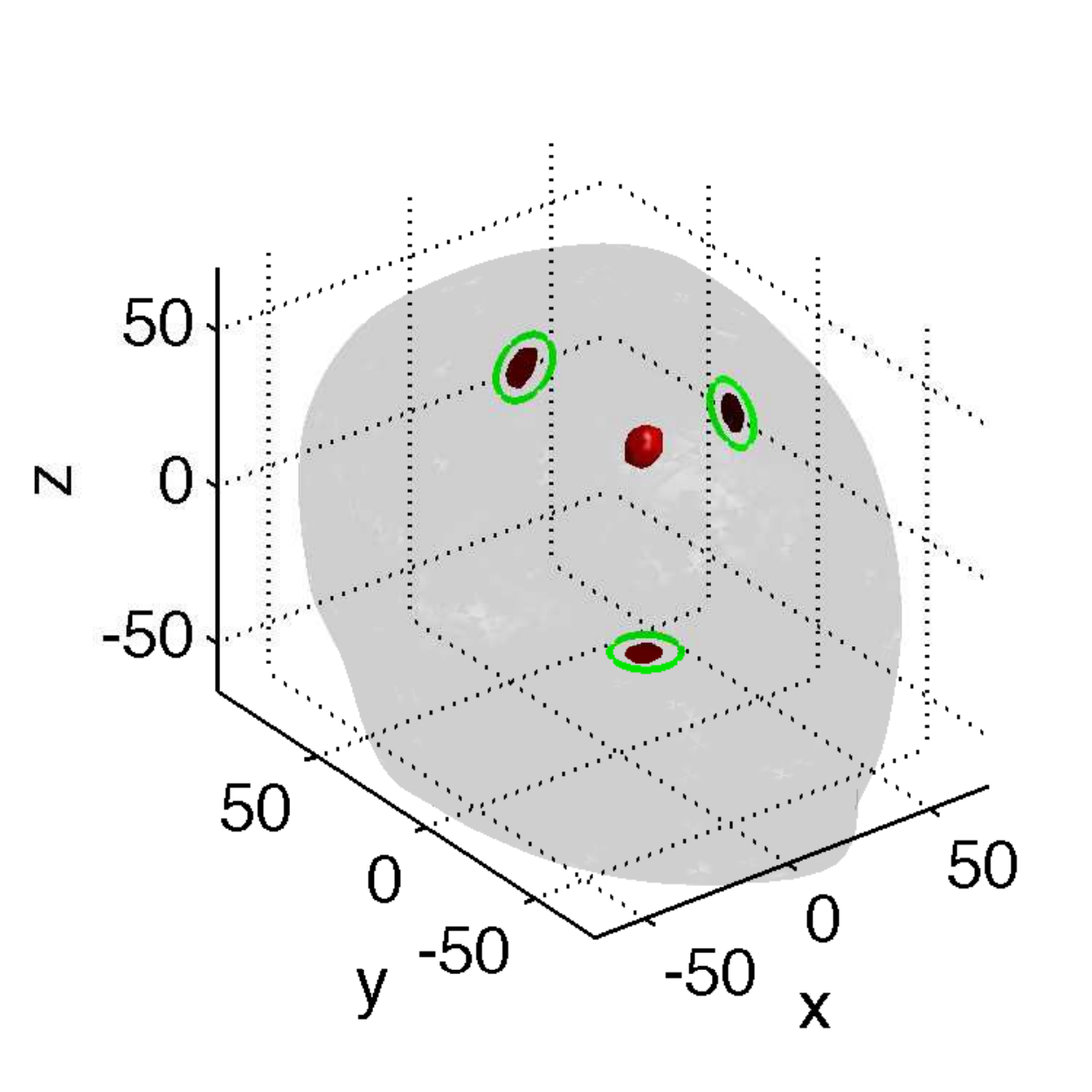}\label{fig:hsback}} \hspace{1cm}
\subfigure[One inclusion in the bask side using $\wt\phi_z^d$; $\wt E_c = 0.13$.]{\includegraphics[width=.34\textwidth]{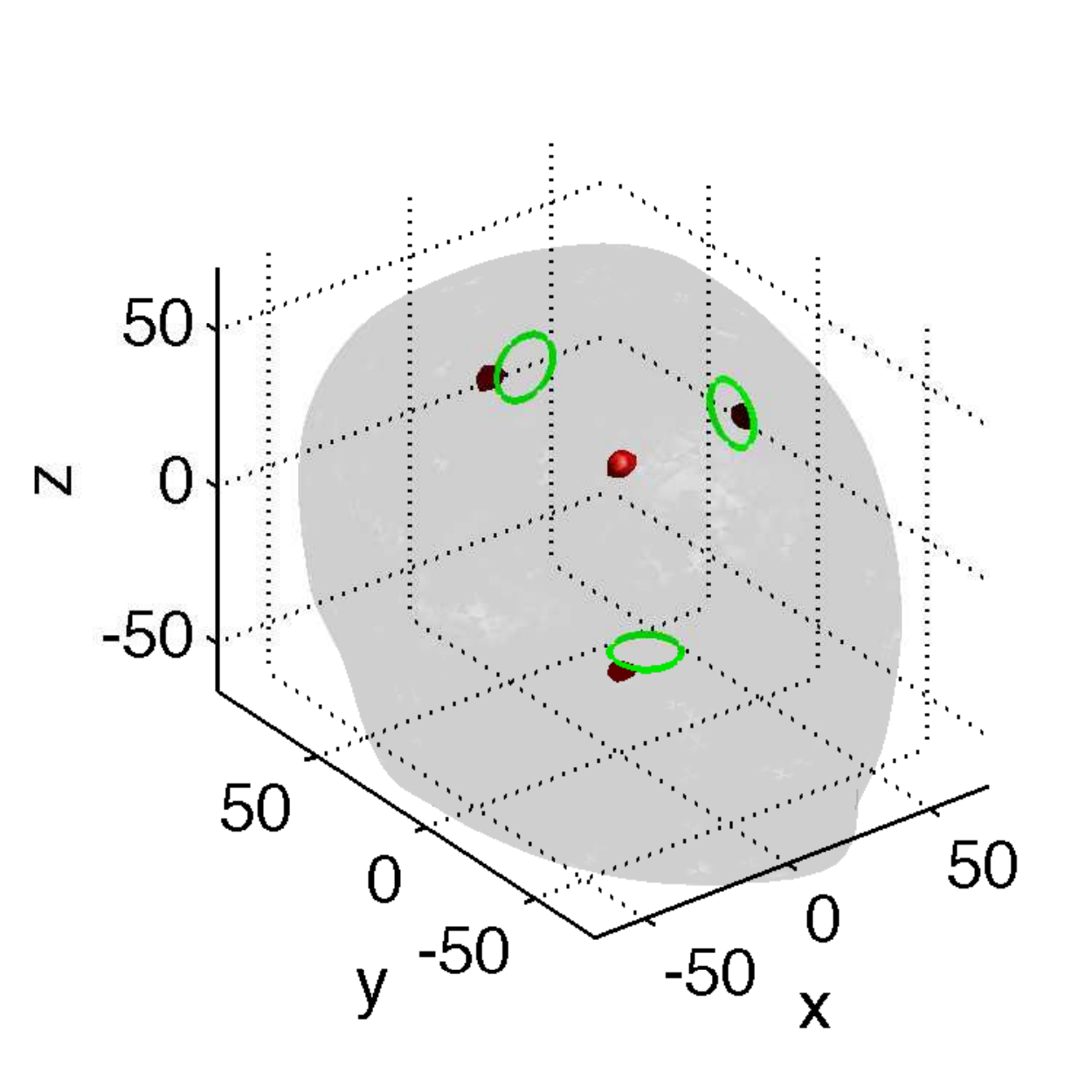}\label{fig:hsbackb}} 
\caption{Reconstruction for a head shape geometry with inhomogeneous background by adding $\delta=1\%$ of noise to the data (same colours as in Figure \ref{fig:cylinder1}).}\label{fig:headinhomogeneous}
\end{figure}

In Figures \ref{fig:cylinder1} to \ref{fig:headinhomogeneous} we 
show reconstructions for various geometries, for constant or piecewise constant 
background conductivities and for inclusions located at different locations.  In all cases, the 
inclusion is a sphere (of radius 1 for the cylinder and 10 for the head shape) and the conductivity value in this inclusion is double that 
 of the background. 
 The different positions of the centre of the sphere are reported in Table \ref{tab:pos}. In order to compare the 
 reconstructed object with the true one, we plot in green the projection of the true object  on the different planes delimiting the plotting region. 
We also choose a quantitative estimate of the quality of the reconstruction by introducing the relative error on the
 location of the barycenter which is defined by
\[
E_c := \frac{|C^{{\rm true}}-C^{{\rm Est}}|}{{\rm diam}(\Om)} \quad {\rm and} \quad \wt E_c := \frac{|C^{{\rm true}}-\wt C^{{\rm Est}}|}{{\rm diam}(\Om)}
\]
where $C^{{\rm true}}$  is the barycenter of $D$, $C^{{\rm Est}}$ (respectively $\wt C^{{\rm Est}}$) is the barycenter of the region 
delimited by the closed iso-surface of level $0.9$ of the normalised indicator function obtained from 
$\rm Ind$ (respectively from $\rm \wt{Ind}$). Finally, ${\rm diam}(\Om)$ stands for the diameter of the computational domain $\Om$. 
 The errors on the location of the barycenter for the different tests presented in this section are reported 
 in the caption of the plots. When  $D$ is not simply connected, we compute $E_c $ and $\wt E_c $
  for each simply connected component and we give the mean of the errors obtained for the different objects. 
 In all this experiments, the noise added to the simulated data is taken such  that $\delta=1\%$ in \eqref{eq:delta}. 
\begin{table}
\centering
\begin{tabular}{|c|c|}
\hline
Figures & Position \\ \hline
Figures \ref{fig:cyl21}, \ref{fig:cyl21b} and  \ref{fig:cyl31}  & (0,5,2) \\
\hline
Figures \ref{fig:cyl22}, \ref{fig:cyl22b} and  \ref{fig:cyl32}  & (0,5,2)  and (5,-2,2) \\
\hline
Figures \ref{fig:hmiddle}, \ref{fig:hsmiddle} and \ref{fig:hsmiddleb}& (0,0,0) \\
\hline
Figures \ref{fig:hback}, \ref{fig:hsback} and \ref{fig:hsbackb}& (40,40,0) \\\hline
\end{tabular}
\caption{Centre of inclusions for the different experiments in section \ref{sec:prelim}.}
\label{tab:pos}
\end{table}
\paragraph*{First setting: cylinder with one ring of electrodes (Figure \ref{fig:cylinder1})\\} We consider a cylindrical domain $\Om$ of radius $10$ and height $7$ with one  ring of electrodes containing 
$32$ electrodes (see Figure \ref{fig:cyl2setup}). The space V of input currents is made of $16$ independent vectors of $\bR^{32}$ corresponding to the so-called 
opposite current pattern (see 
\cite[Chapter 12]{MueSil12} for more details on different current patterns). This means that all the electrodes are set to $0$ except for two of them that are 
geometrically opposite to each other. One of these is set to $1$, the other one to $-1$ so that the constraint on the input current is satisfied. In this experiment, the 
background conductivity is taken constant and equal to $1$ while the contact impedance value is $z=5$ for all the electrodes (we keep this value for
all experiments).  In Figure \ref{fig:cyl21} we see that with one inclusion, the 
algorithm with $\phi_z^d$ as test function finds the correct (x,y) location of the inhomogeneity but not the $z$ location whereas the algorithm 
with $\tilde \phi_z^d$ (Figure \ref{fig:cyl21b})
does not even find the $(x,y)$ location. The reason why Ind gives the correct $(x,y)$ location
 but not the $z$ one is because we only have one ring of electrodes at  the same $z$ position. 
We also run an experiment with two well separated inclusions and the algorithm with $\phi_z^d$ (Figure 
\ref{fig:cyl22}) seems to give 
quite accurate results in this case as well while the use of $\tilde \phi_z^d$ (Figure \ref{fig:cyl22b}) only gives a rough idea of the location
of the objects. Therefore, in the following we will mainly use the indicator function Ind.

\paragraph*{Second setting: cylinder with two rings of electrodes (Figure \ref{fig:cylinder2})\\} 
This experiment is similar to the first one (the domain $\Om$ is the same) except  that we have two rings of $20$ 
electrodes each. The injection protocol is again a pairwise injection  which corresponds to the opposite injection pattern for each ring of electrodes. Then, the data consist 
of a $40\times20$ matrix. The results are similar to the previous case, except that the resolution in the $z$ 
direction in much better in this case since we accurately find the $(x,y,z)$ location of the inhomogeneity if we use the correct dipole function (Figures \ref{fig:cyl31} 
and \ref{fig:cyl32})
This experiment illustrates that it is not necessary to have electrodes covering the entire domain 
to have good quality reconstructions since we do 
not have any electrode on the top and on the bottom of the cylinder and still we find the $z$ location with a good accuracy. 

\paragraph*{Third setting: human head with homogeneous background (Figures \ref{fig:headhomogeneous}) \\} 
We use a head shape domain with $31$ electrodes that covers  a portion of the physically 
accessible part of the head (see Figure \ref{fig:headhmosetup}). In this setting,  we show that the Factorization method still performs
well.
The protocol injection is very similar to the previous one and the data correspond to a 
$31\times20$ 
Neumann-to-Dirichlet matrix. We tried two positions for the inclusion which is a sphere of radius $10$ and of conductivity $2\sigma_0$ with a background 
conductivity of $\sigma_0=1$. We still use a constant contact impedance which is $z=5$ for all the electrodes.
 In both cases the location of the inclusion is accurately found. 

\paragraph*{Fourth setting: human head with inhomogeneous background (Figure \ref{fig:headinhomogeneous})\\}
We use the same head shape geometry as previously, but this time the conductivity $
\sigma_0$ of the background is not constant.  
We take a piecewise constant conductivity with values $1.5 \times 10^{-4}$ in the yellow part, $2.0 \times10^{-5}$ in 
the green part and $ 4.4\times10^{-4}$ in the red part of the domain (see Figure \ref{fig:skullsetting}).
 These values correspond approximately to the conductivity of 
the skin, the skull and the brain respectively of a human head (see \cite{Hor06} and references therein).
The inclusion is still a sphere of radius $10$ and of conductivity $2\sigma_0=8.8\times 10^{-4}$ located at the 
 same places as in the previous setting.  Let us mention the fact 
that for the two locations, the inclusion is inside the brain (the red part of the computational domain). 
The reconstructions with the exact dipole tests functions are very precise but when we use  
the dipole function of the free space we observe a misplacement of the reconstructed object. 
This difference is actually significant when the inclusion is close the inhomogeneous layers (compare Figure \ref{fig:hsback} with 
Figure \ref{fig:hsbackb}). 
The result obtained with  the free space dipole function  will be used for comparison in section \ref{sec:model}.

\subsection{Influence of the truncation value and of the inclusion's size on the reconstruction}
\begin{figure}
\centering
\subfigure[Iso-surface value: $0.6$ ; $E_c = 0.02$.]{
\includegraphics[width=.31\textwidth]{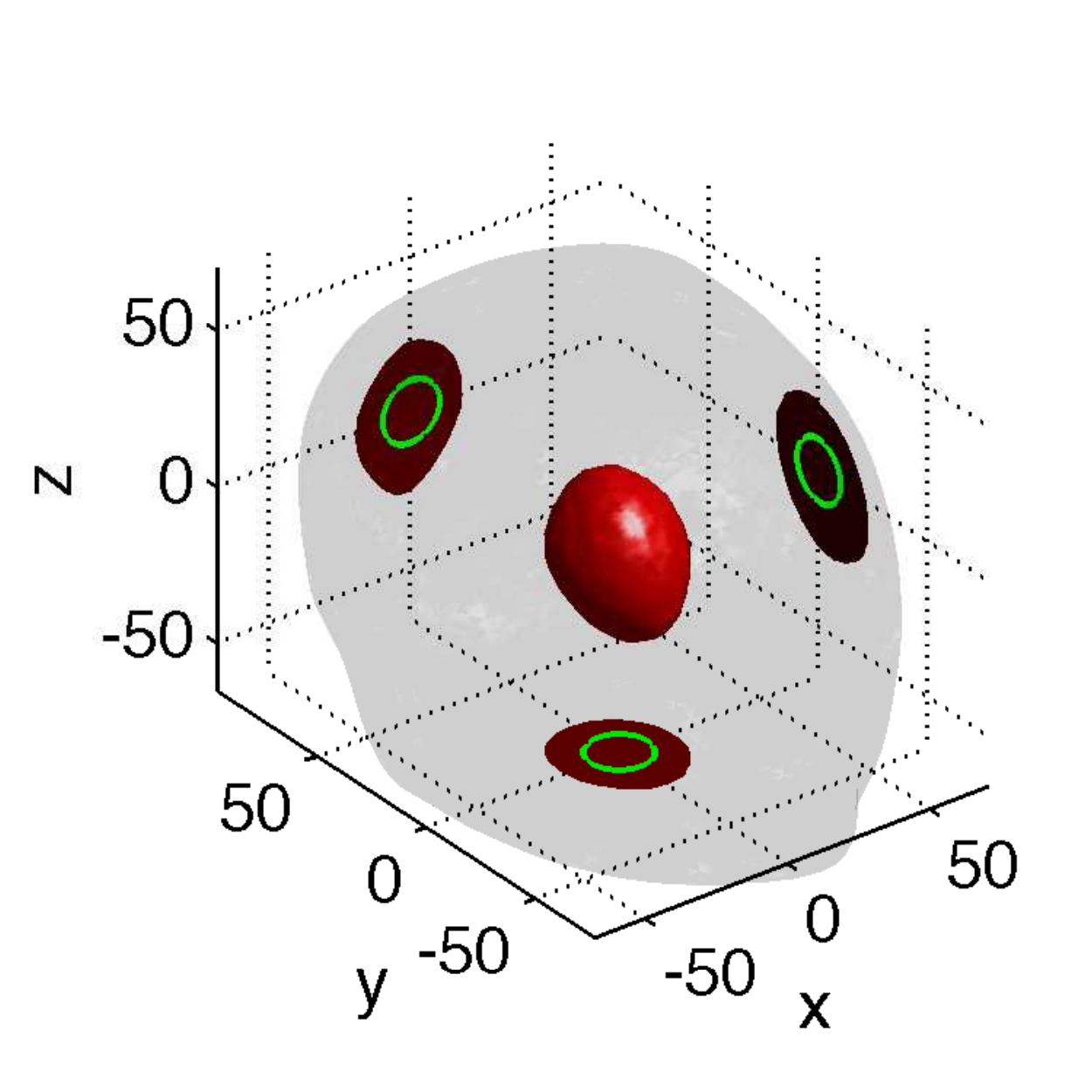}}\hfill
\subfigure[Iso-surface value: $0.8$ ; $E_c = 0.01$.]{
\includegraphics[width=.31\textwidth]{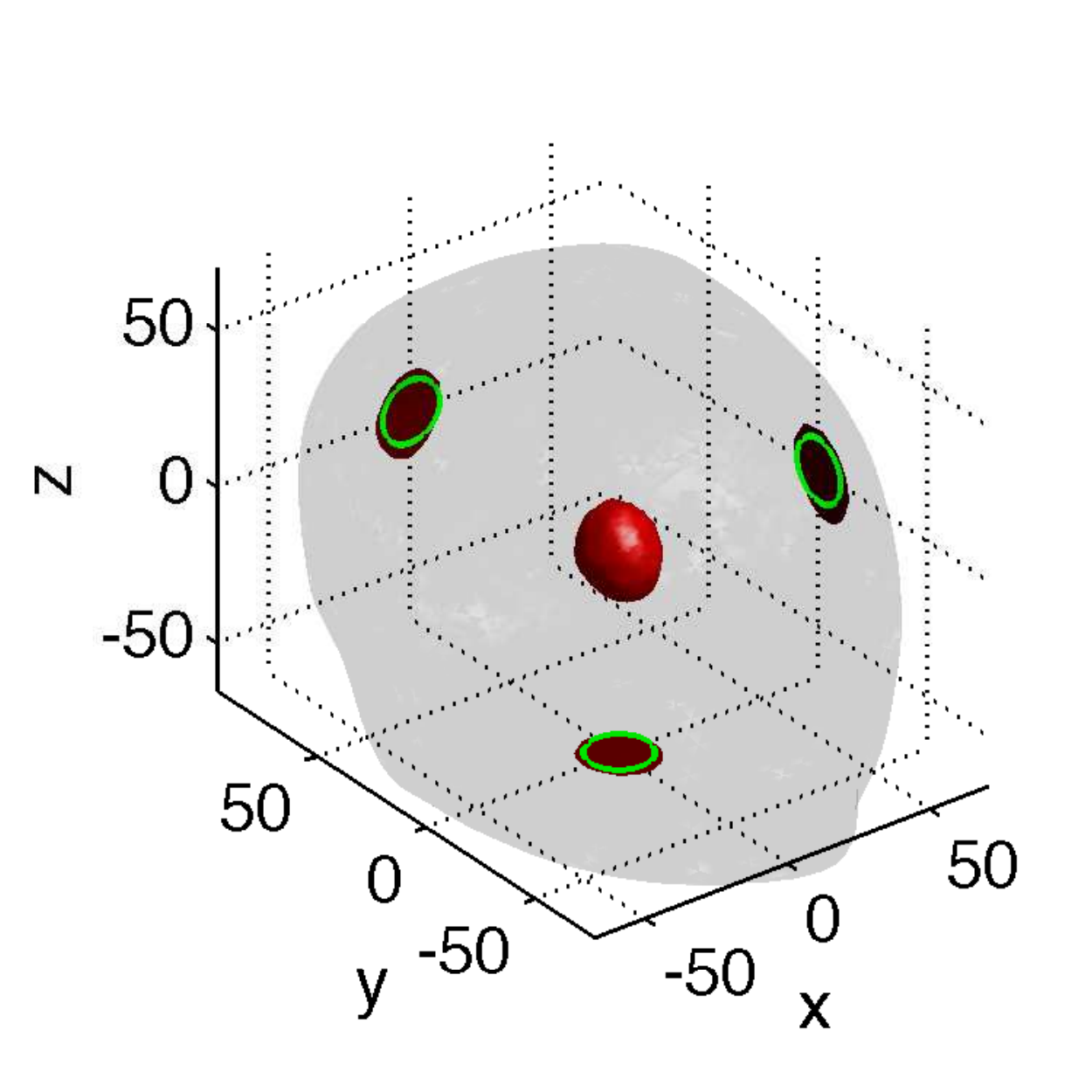}} \hfill
\subfigure[Iso-surface value: $0.95$ ; $E_c = 0.01$.]{
\includegraphics[width=.31\textwidth]{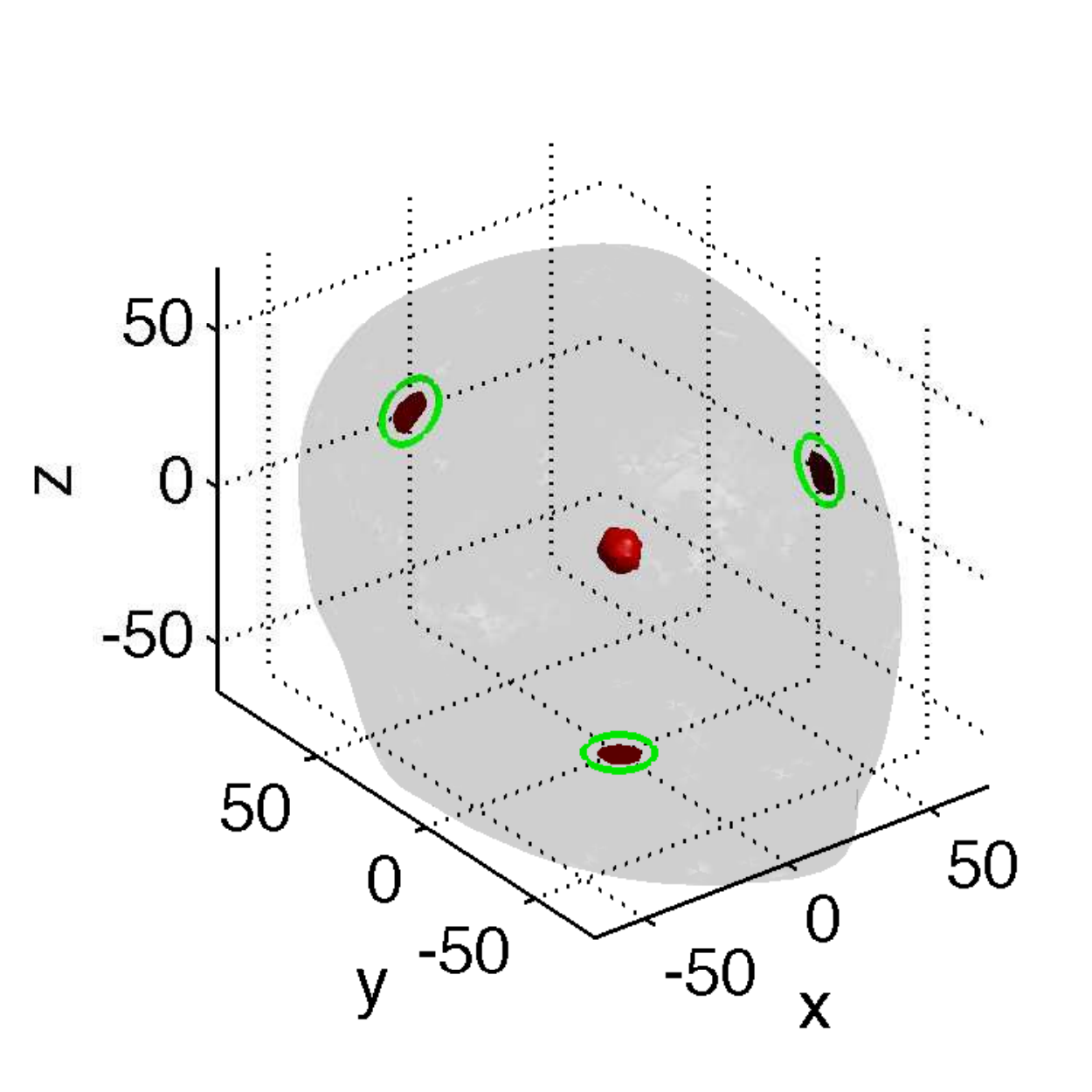}} 
\caption{Reconstruction of a centred inclusion by using $\phi_z^d$ for different choices of iso-surface  of function 
Ind (same colours as in Figure \ref{fig:cylinder1}).}\label{fig:headtruncation}
\end{figure}

\begin{figure}
\centering
\subfigure[Radius of inclusion: $5$; $E_c = 0.01$.]{
\includegraphics[width=.31\textwidth]{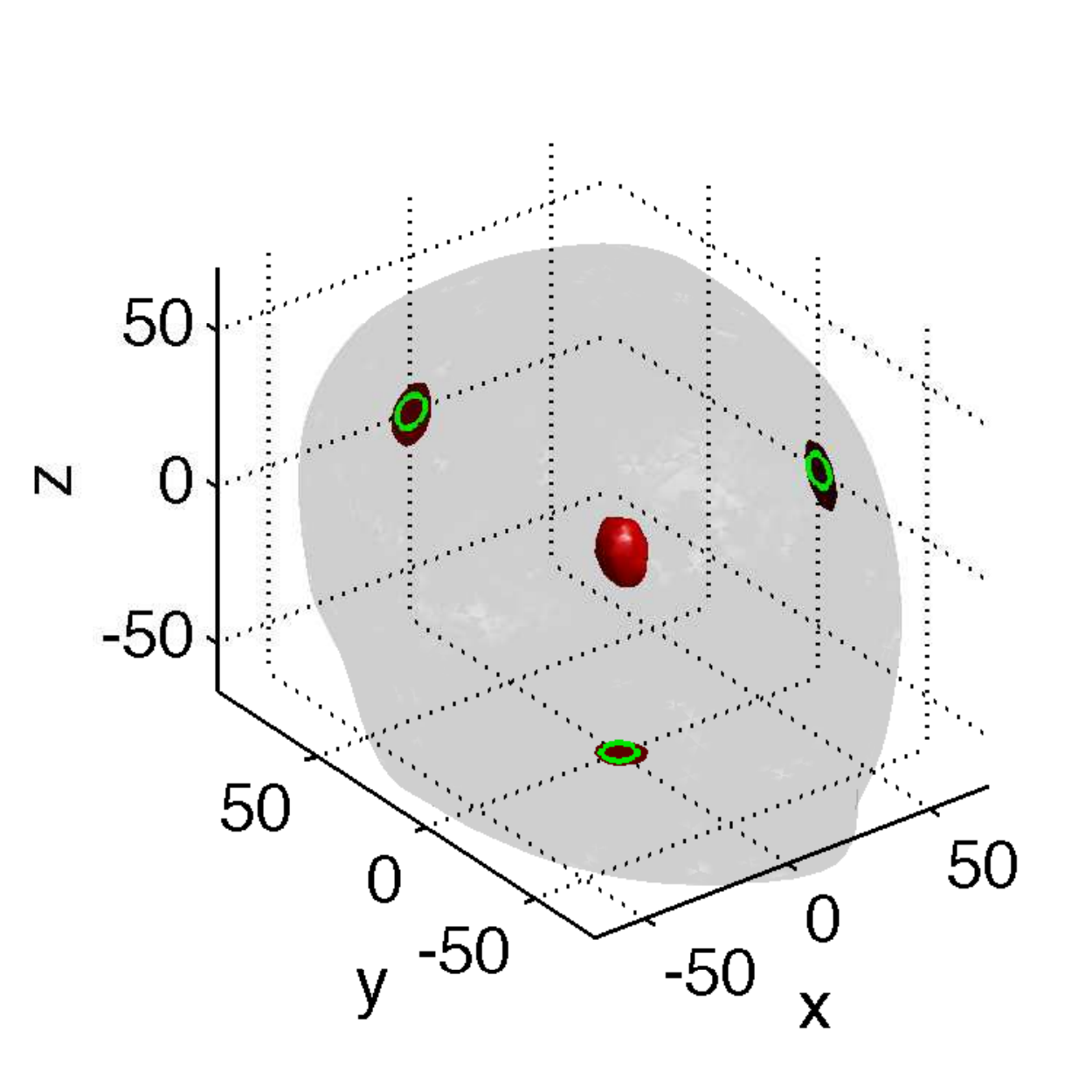} \label{fig:small}}\hfill
\subfigure[Radius of inclusion: $15$; $E_c = 0.01$.]{
\includegraphics[width=.31\textwidth]{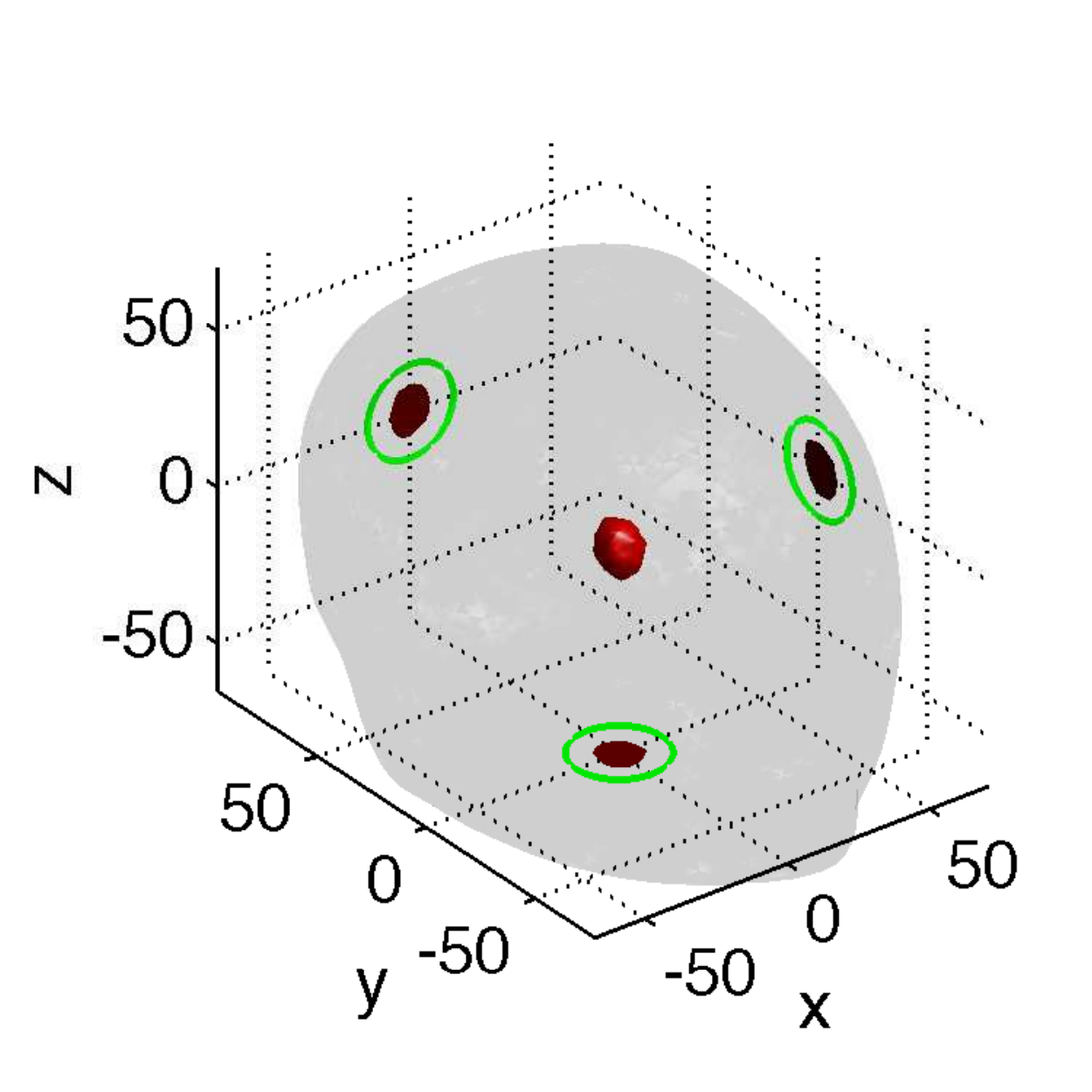}} \hfill
\subfigure[Radius of inclusion: $20$; $E_c = 0.005$.]{
\includegraphics[width=.31\textwidth]{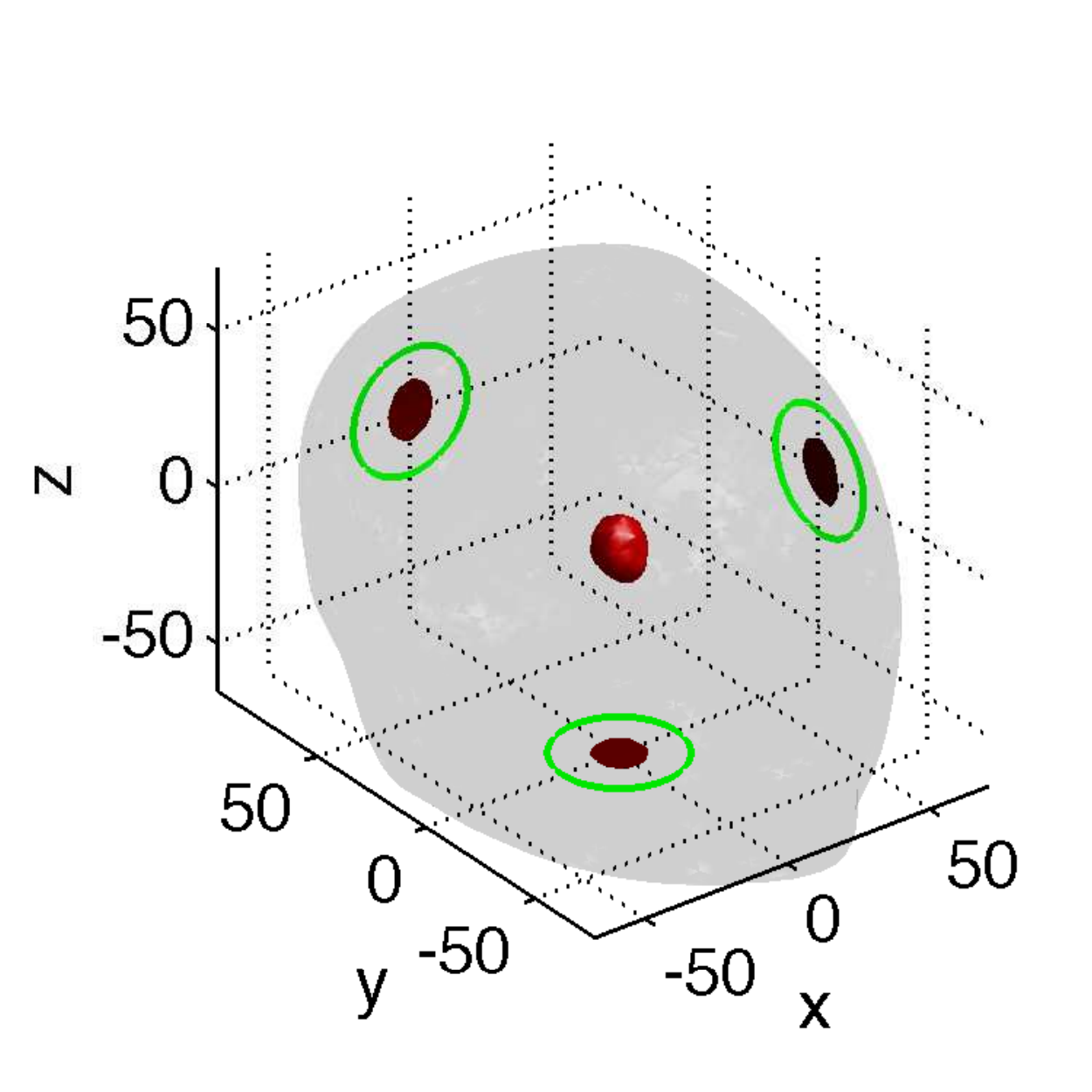}} 
\caption{Reconstruction of a centred inclusion of different sizes by using $\phi_z^d$
 (same colours as in Figure \ref{fig:cylinder1}).}\label{fig:headsize}
\end{figure}

In the next set of experiments we illustrate the influence of the choice of iso-surface (Figure \ref{fig:headtruncation}) and of
the size of the inclusion (Figure \ref{fig:headsize}). All these simulations are performed in the same geometry as in Figure \ref{fig:headinhomogeneous}
and the conductivity value of the background and of the inclusion are also the same. Therefore, one can compare these results to the 
reference results presented in Figure \ref{fig:headinhomogeneous}.

Regarding the influence of the value of the iso-surface that defines the boundary of the inclusion, we remark that it greatly affects 
the size of the reconstructed object but not its location. As it has already been observed in previous work, 
in its actually state, the Factorization  method is probably not the right tool
to estimate the size of a defect since this quantity depends too much on the choice of the iso-surface value. Nevertheless, whatever the
 value of the iso-surface, the object we reconstruct has the correct shape and the correct location.

These remarks also apply to our second test (Figure \ref{fig:headsize}) where we plot the iso-surface of value 0.9
but this time for inclusions of different sizes. Nevertheless, let us mention the fact that even for a small inclusion (Figure \ref{fig:small}) 
the method locates accurately the defect.

\subsection{Robustness to modelling errors}
\label{sec:model}
We conclude our numerical analysis of the Factorization method in the context of EIT for brain imaging 
by studying the influence of modelling errors and noise on the reconstructions. 

\begin{figure}
\centering
\subfigure[Noise level: $\delta=3\%$; $E_c = 0.02$.]{
\includegraphics[width=.31\textwidth]{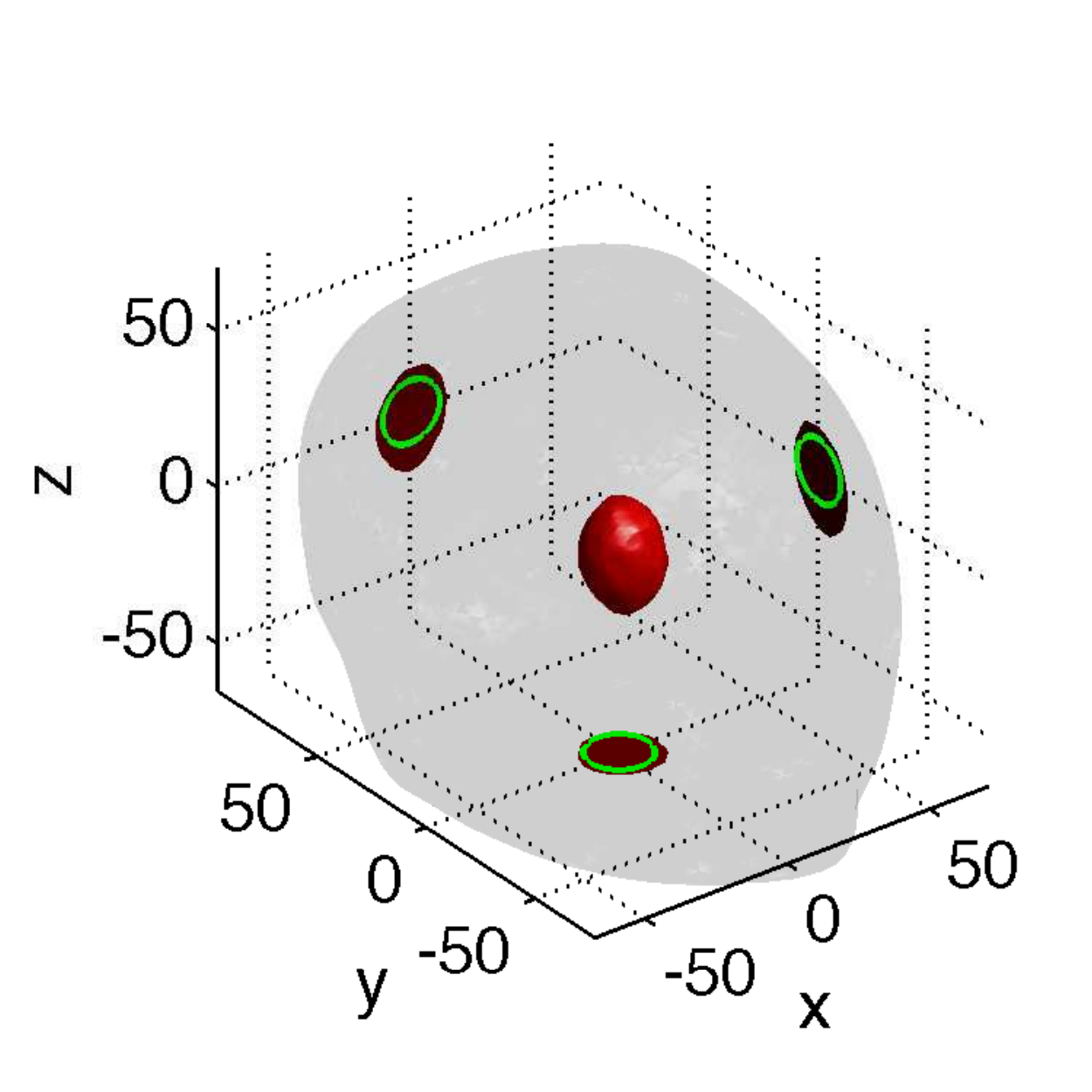}}\hfill
\subfigure[ Noise level: $\delta=7\%$; $E_c = 0.04$.]{
\includegraphics[width=.31\textwidth]{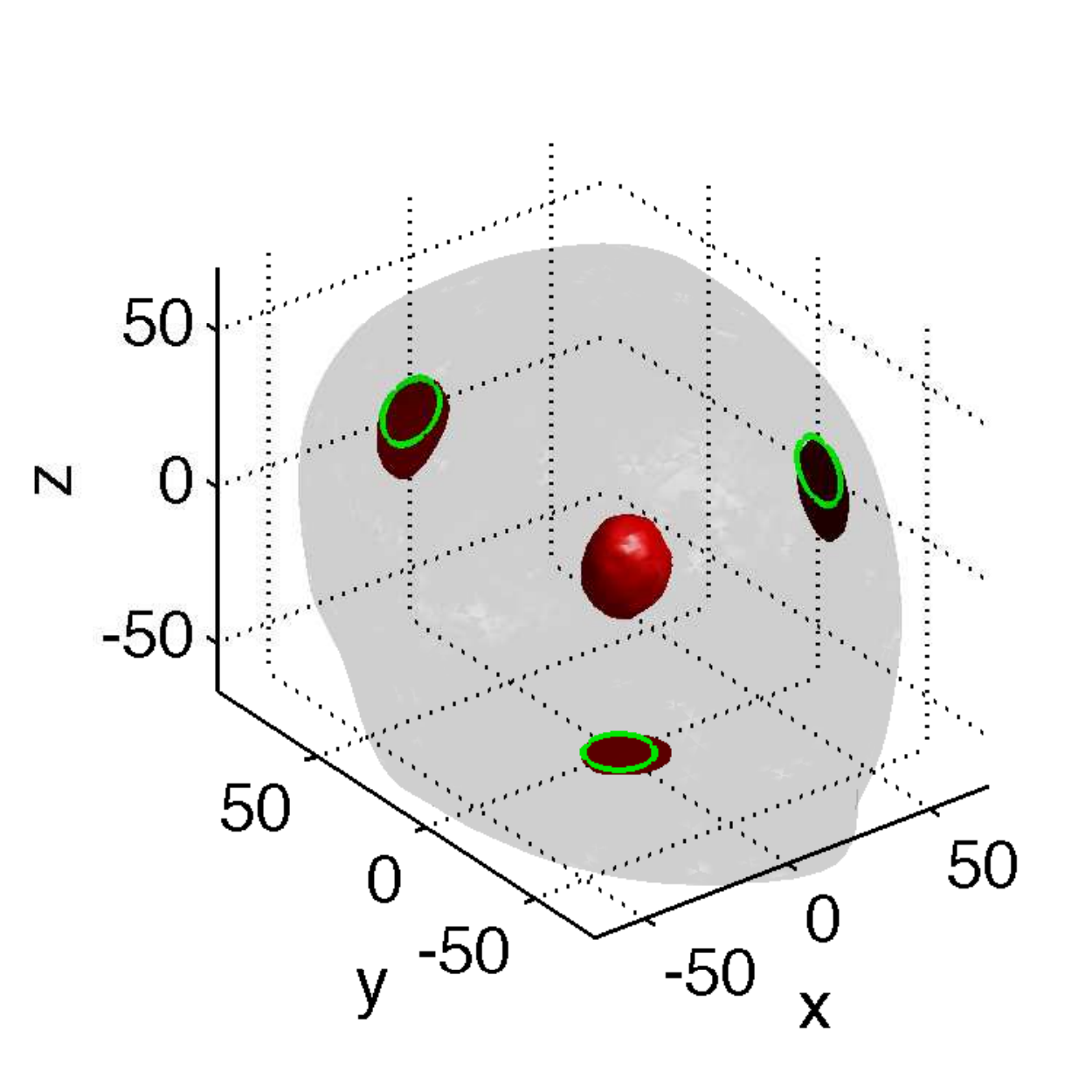}} \hfill
\subfigure[Noise level: $\delta=10\%$; $E_c = 0.05$.]{
\includegraphics[width=.31\textwidth]{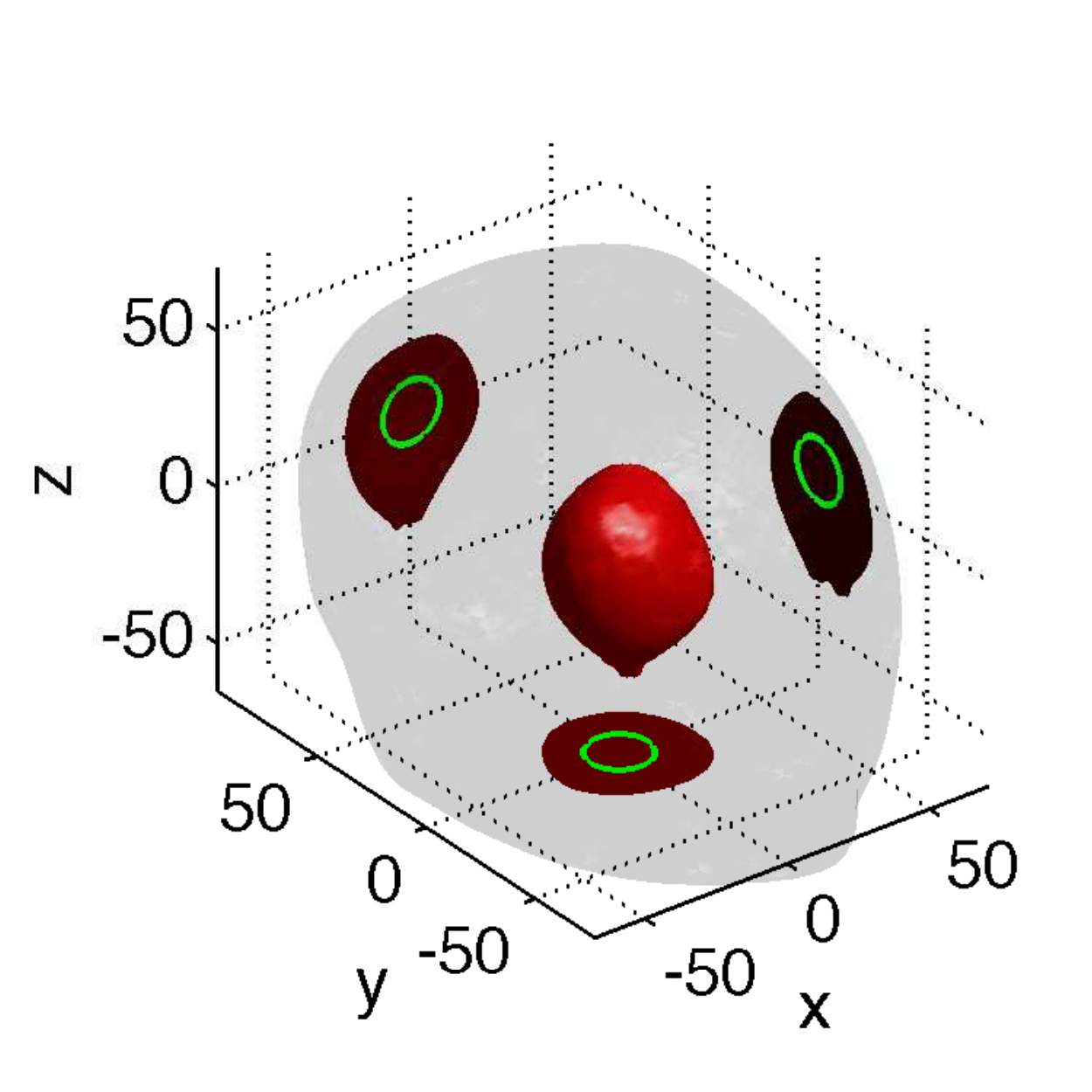}\label{fig:noisehigh}} 
\caption{Reconstruction of a centred inclusion by using $\phi_z^d$ with various level of noise
 (same colours as in Figure \ref{fig:cylinder1}).}\label{fig:headnoise}
\end{figure}

First of all, in Figure \ref{fig:headnoise} we repeat the same experiment as in Figure \ref{fig:hsmiddle} and we plot
the iso-surface of value 0.9 of the indicator function Ind for different level of noise added to the simulated data.
The size of the reconstructed defect is strongly affected by the noise level but even for $10\%$ of noise we
still obtain a very accurate estimate of the location of the inclusion (see Figure \ref{fig:noisehigh}).

\begin{figure}
\centering
\subfigure[Slices of $\Om$ (in blue) and of $\Om_\ve$ (in orange).]{
\raisebox{1cm}{\includegraphics[width=.22\textwidth]{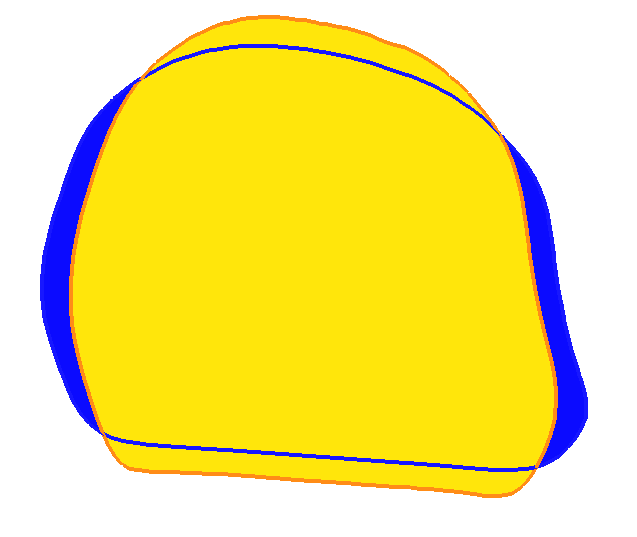}}\label{fig:headpertshape}}\hfill
\subfigure[ One inclusion in the middle using $\phi_z^d$; $E_c = 0.05$.]{
\includegraphics[width=.34\textwidth]{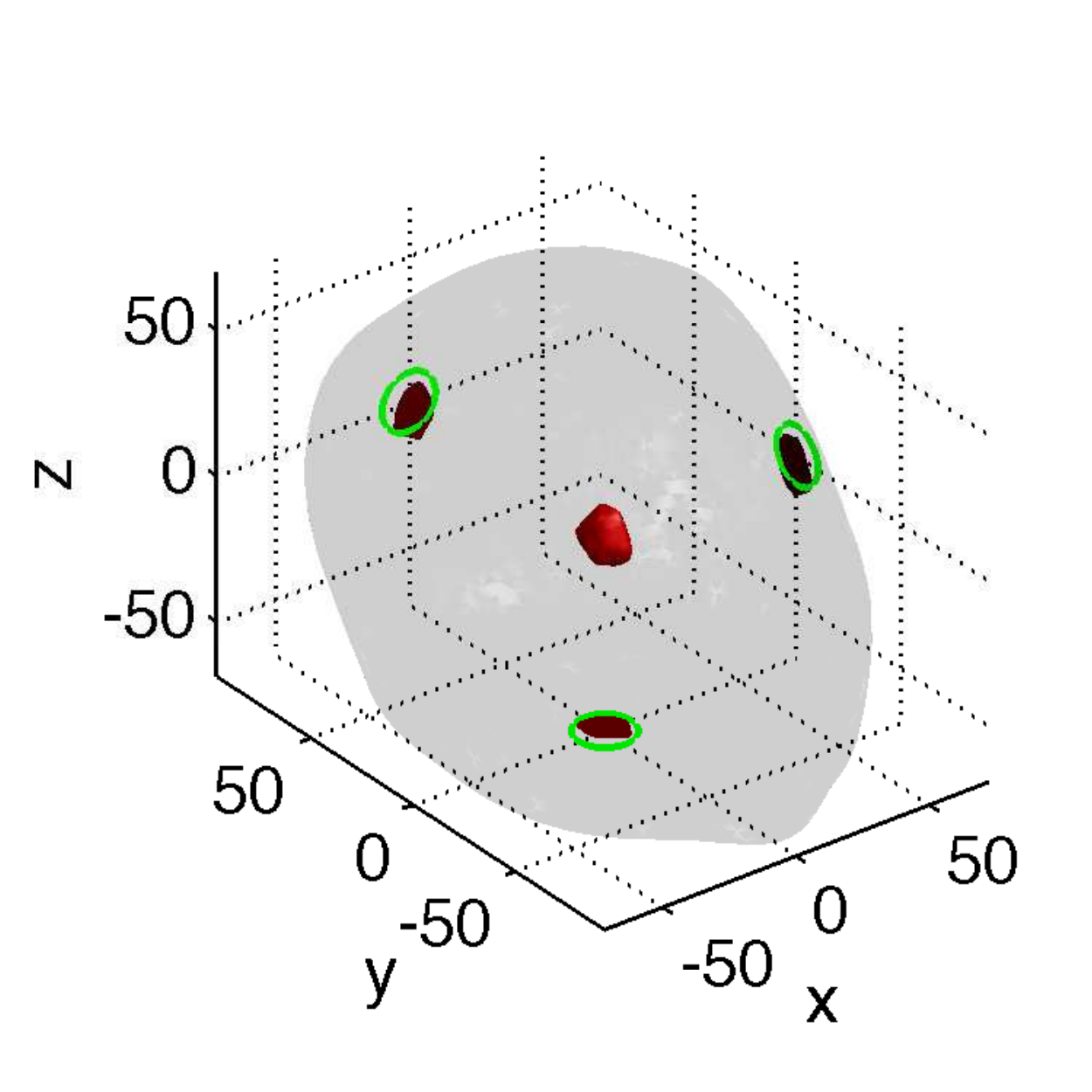}} \hfill
\subfigure[One inclusion in the  back side using $\phi_z^d$; $E_c = 0.08$.]{
\includegraphics[width=.34\textwidth]{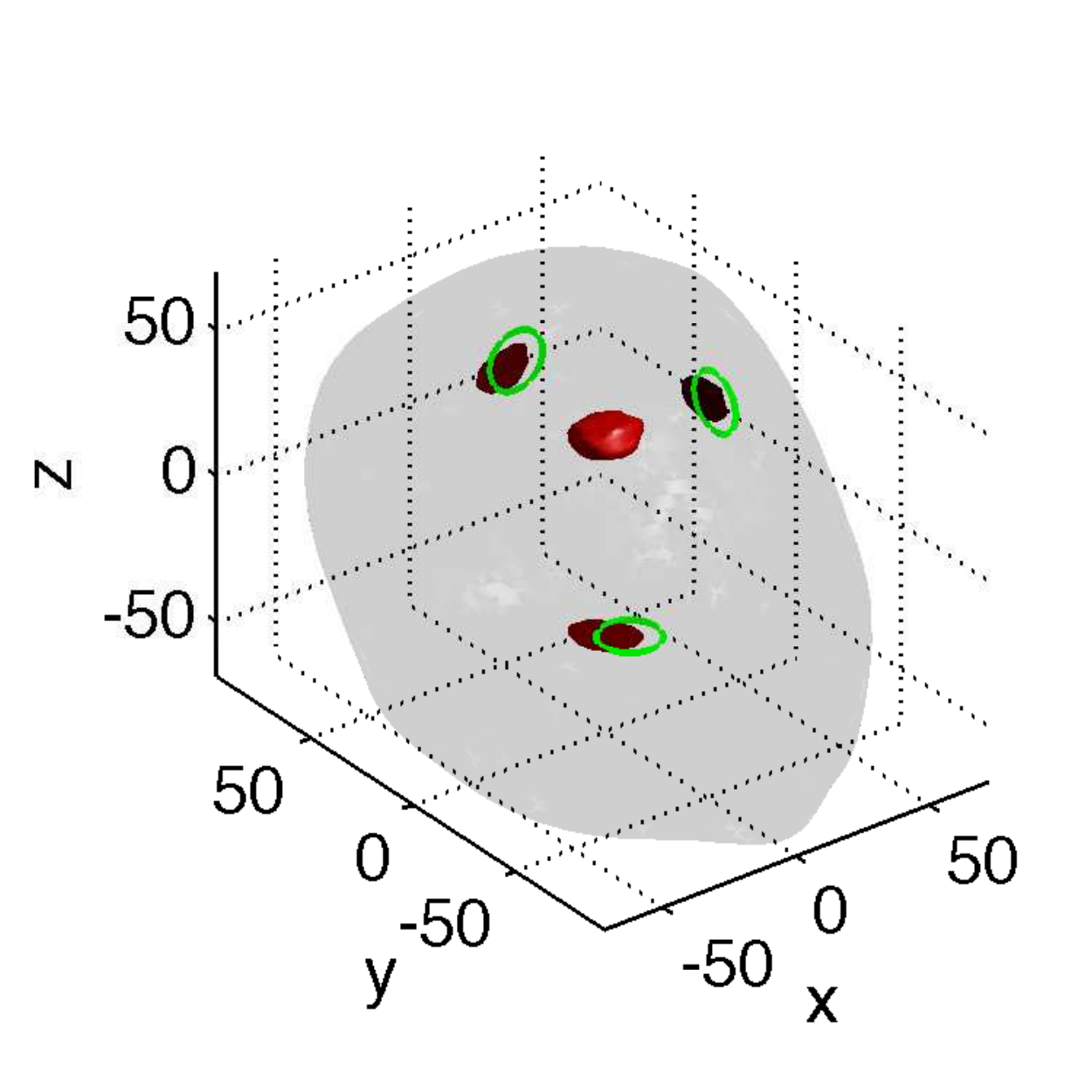}} 
\caption{Reconstruction for a deformed head shape geometry with homogeneous background 
by adding $\delta=1\%$ of noise to the data (same colours as in Figure \ref{fig:cylinder1}).}\label{fig:headperturbed}
\end{figure}

We also test the sensitivity of the method with respect to errors on the shape of the domain $\Om$ and on the electrode's placement. Indeed, this is of 
crucial importance for applications since one usually has only a rough idea of the shape of $\Om$ and of the electrode positions. 
We consider in this experiment the homogeneous head that has already been used in Figure \ref{fig:headhomogeneous}. We choose to not
use the inhomogeneous head of Figure \ref{fig:headinhomogeneous} since we would like 
to decouple errors in the background conductivity and errors in the shape of the domain.
To simulate such modelling errors, we deform the
boundary of the domain $\Om$ by expanding it in the direction $x$ and $z$ and by contracting it in the direction $y$ and to obtain a perturbed domain 
$\Om_\ve$. See Figure \ref{fig:headpertshape} for a superposition of  cut views of $\Om$ and $\Om_\ve$, in blue we plot a slice of $\Om$ and in orange 
a slice of $\Om_\ve$. We remark that by doing so we also generate a domain with
inexact electrode positioning. Let us denote $E_\ve \in \bR^{31\times3}$, respectively $E \in \bR^{31\times3}$,
 the center of the electrodes corresponding to $\Om_\ve$, respectively $\Om$.
 We then simulate the Neumann-to-Dirichlet matrix  data in the domain 
$\Om_\ve$ with inexact electrode positions $E_\ve$
 and use the Green's function of $\Om$ evaluated at the exact electrode position $E$ to produce an image. In Figure \ref{fig:headperturbed}
  we report results  for $\Om_\ve$ and $E_\ve$ being such that the mean relative  errors on the shape and the electrode positions are 
  \[
  \frac{\int_{\partial \Om} {\rm dist}(x,\partial \Om_\ve) dx}{{\rm diam}(\Om)} = 7\% \quad {\rm and} 
  \quad \frac{\sum_{i=1}^{31} |E^i-E^i_\ve |}{{\rm diam}(\Om)} = 11\%
  \]
  where ${\rm dist}$ stands for the distance function and $E^i \in \bR^3$, respectively $E^i_\ve \in \bR^3$, 
  contains the position of the $i^{th}$ electrode on $\partial \Om$, respectively on $\partial \Om_\ve$. 
  The main conclusion we can draw from these experiments is that the Factorization method we presented 
  in this paper is stable with respect to errors on the shape of the measurement domain.

\begin{figure}
\centering
\subfigure[$\gamma = 0.1$ in \eqref{eq:sigmaeta}; $E_c = 0.03$.]{
\includegraphics[width=.31\textwidth]{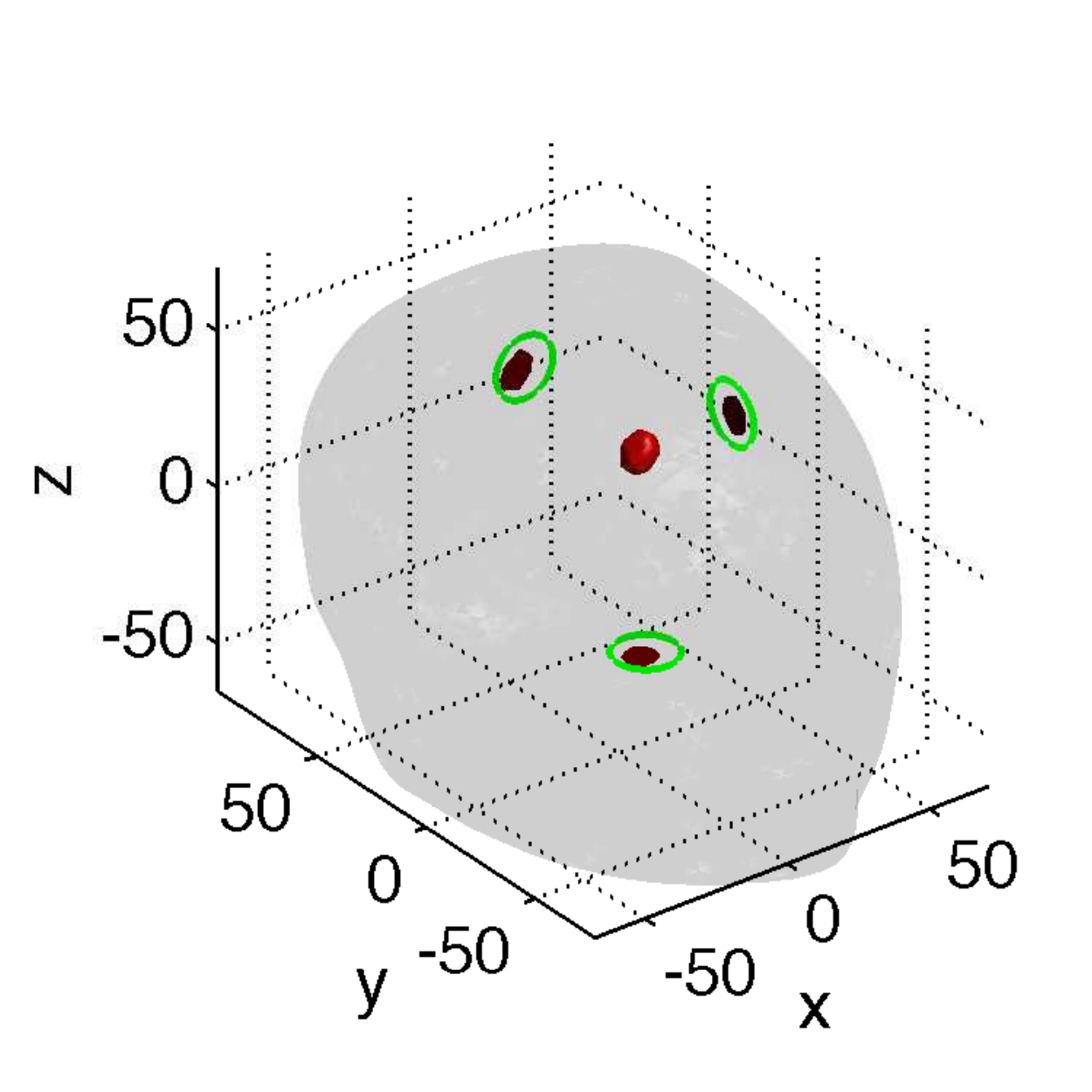}}\hfill
\subfigure[$\gamma = 0.2$ in \eqref{eq:sigmaeta}; $E_c = 0.04$.]{
\includegraphics[width=.31\textwidth]{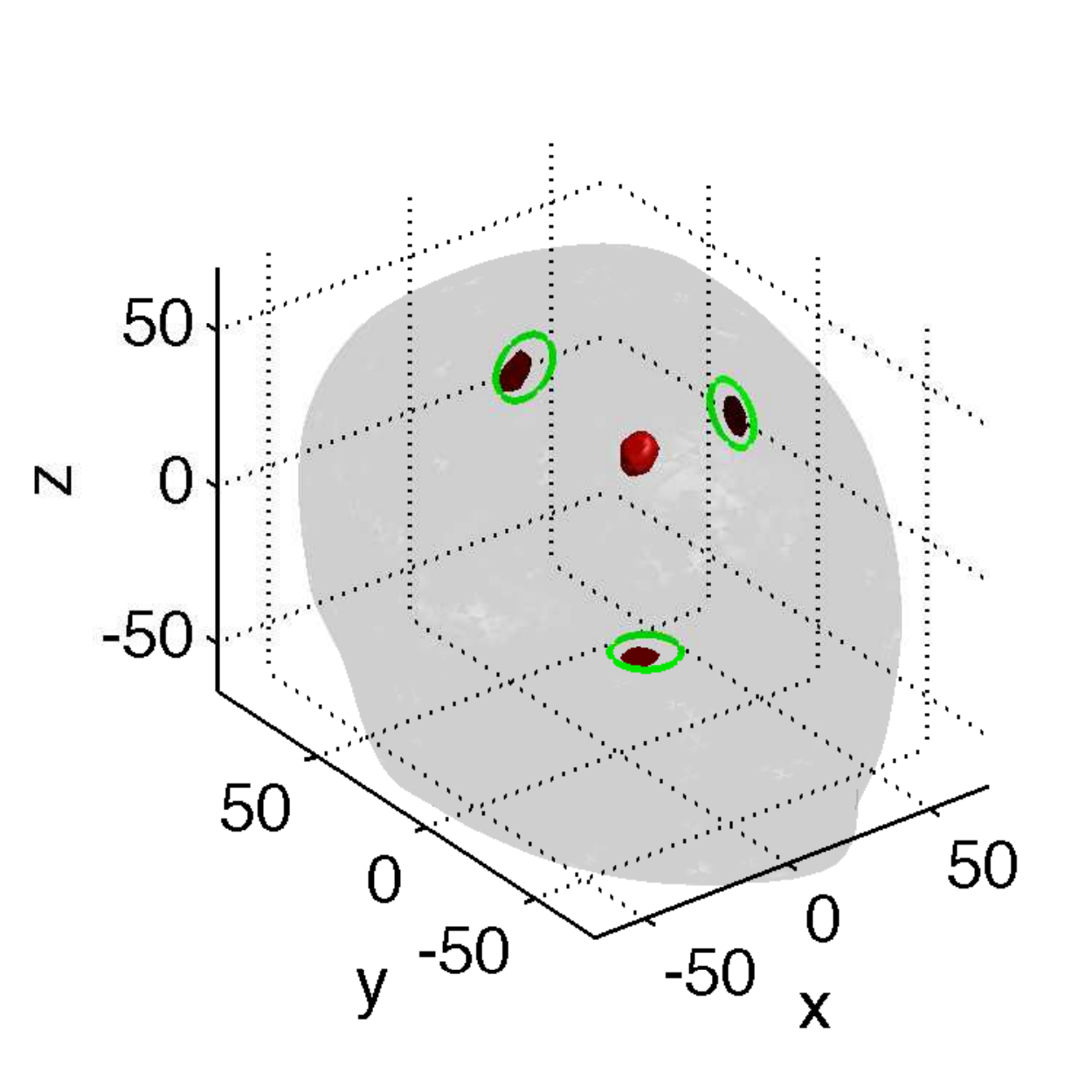}} \hfill
\subfigure[$\gamma = 0.5$ in \eqref{eq:sigmaeta}; $E_c = 0.09$.]{
\includegraphics[width=.31\textwidth]{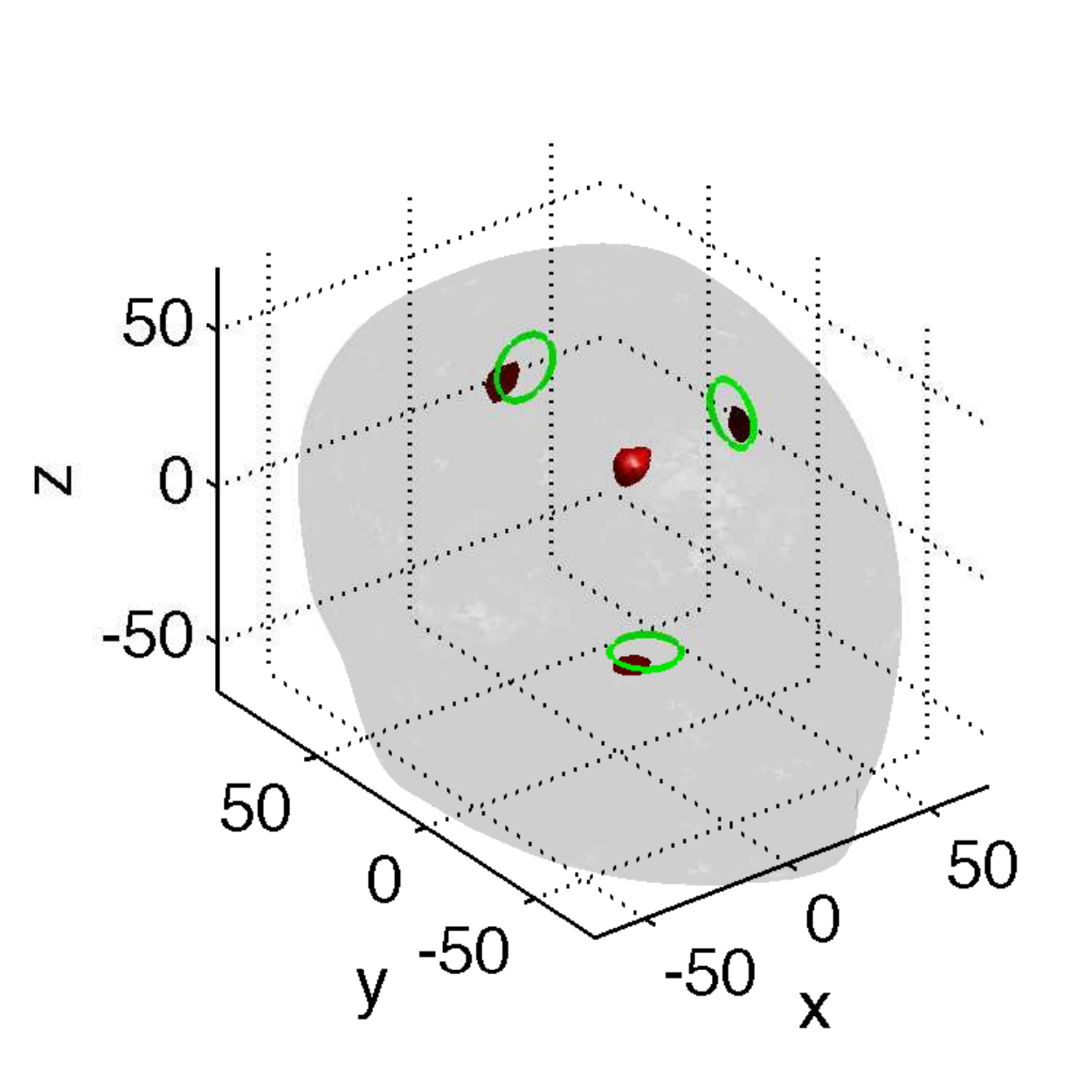}} \label{fig:hbackb}
\caption{Reconstruction of a lateral inclusion by using $\phi_z^d$ with error on the background conductivity
 (same colours as in Figure \ref{fig:cylinder1}).}\label{fig:headsigma}
\end{figure}

The last source of modelling errors that certainly affects the Factorization method is in the background conductivity value.
For this last experiment we go back to the inhomogeneous head introduced in Figure \ref{fig:headinhomogeneous}.
We generate data for a noisy background conductivity $\sigma_0^\gamma$ given for $\gamma>0$ by
\begin{equation}
\label{eq:sigmaeta}
\sigma_0^\gamma =\begin{cases} (1+\gamma) \sigma_0 \quad \mbox{in the yellow part,}\\
(1-\gamma) \sigma_0 \quad \mbox{in the green part,} \\
(1+\gamma) \sigma_0 \quad \mbox{in the blue part} 
\end{cases}
\end{equation}
for $\sigma_0$ being as in Figure \ref{fig:headinhomogeneous} and the colours refer to the regions in Figure \ref{fig:skullsetting}. 
Then, we build the indicator 
function Ind by using the dipole function associated with the background conductivity $\sigma_0$. 
In Figure \ref{fig:headsigma} we plot the iso-surface of value 0.9 of the indicator function Ind for an inclusion located at $(40,40,0)$.
In view of  Figure \ref{fig:headinhomogeneous} we expect that this position is  more affected by error
in   the background conductivity than the middle position since in the case of the middle position (Figures \ref{fig:hsmiddle} and \ref{fig:hsmiddleb})
even the indicator function $\wt {\rm Ind}$ gave accurate reconstructions.
The results reported in Figure \ref{fig:headsigma} show that up to $20\%$ of noise on the background conductivity, the Factorization method 
still performs reasonably well and gives an accurate estimate of the object location. Let us also keep in 
mind that in situations where one does not have any knowledge of the background conductivity, the free space dipole functions are 
in practise good candidates to probe the domain (see Figures \ref{fig:hsmiddleb}and \ref{fig:hsbackb}).

\section{Conclusions}
The main focus of this paper was to show via numerical experiments that  the Factorization method 
can be used to solve realistic three dimensional imaging problem with EIT data. 
We presented a possible implementation of the computation of the Neumann Green's 
function for three dimensional problems and we provided   reconstructions  obtained with 
the Factorization algorithm by using the numerically computed Green's function. In the considered
 test cases the obtained results are rather accurate. The Factorization method  actually finds the 
 location of an inclusion with $1\%$ accuracy in the challenging and realistic case of a head shape 
 domain with homogeneous and inhomogeneous background conductivities from noisy simulated 
 data obtained with $31$ electrodes that cover only a part of the domain's boundary. We have also
 show that the Factorization method is robust with respect to noise on the measurements 
 and various modelling errors such as errors on electrode positioning, on the shape of the domain
and  on the background conductivity.

A direct extension of this work would be to perform an experimental study and to incorporate 
some systematic criterion to choose the truncation level for valuable 
singular values and to pick an iso-surface of the indicator function. We think for example of 
the criterions introduced in \cite{BruHan03}. Another possible extension, which is probably the
 main interest of the Factorization method for applications, is to couple the result obtained by 
 the Factorization algorithm with optimisation algorithms as it is proposed for example in \cite{ChHaSepre}.

\section*{Acknowledgement}
The research of N. C. is supported by the Medical Research Council Grant MR/K00767X/1. The authors are grateful to the EIT group in the 
Department of Medical Physics and Bioengineering of University College London for providing the meshes of the head.
The authors would like to thank  the anonymous referees for their valuable comments and suggestions.

\bibliographystyle{plain}
\bibliography{factorization_EIT}{}

\end{document}